\documentclass[11pt]{article}
	
	\newcommand{\blind}{0}
	
    \makeatletter
    \renewcommand\section{\@startsection {section}{1}{\z@}%
                                       {-3.5ex \@plus -1ex \@minus -.2ex}%
                                       {2.3ex \@plus.2ex}%
                                       {\normalfont\fontfamily{phv}\fontsize{16}{19}\bfseries}}
    \renewcommand\subsection{\@startsection{subsection}{2}{\z@}%
                                         {-3.25ex\@plus -1ex \@minus -.2ex}%
                                         {1.5ex \@plus .2ex}%
                                         {\normalfont\fontfamily{phv}\fontsize{14}{17}\bfseries}}
    \renewcommand\subsubsection{\@startsection{subsubsection}{3}{\z@}%
                                        {-3.25ex\@plus -1ex \@minus -.2ex}%
                                         {1.5ex \@plus .2ex}%
                                         {\normalfont\normalsize\fontfamily{phv}\fontsize{14}{17}\selectfont}}
    \makeatother
	
	\usepackage{amsmath}
	\usepackage{graphicx}
	\usepackage{enumerate}
	\usepackage{natbib} 
	\usepackage{url} 
        \usepackage{xcolor}

        \usepackage{subfigure}
        \usepackage{tikz}
        \usetikzlibrary{arrows.meta,positioning,decorations.pathmorphing}

        \usepackage{booktabs}
        \usepackage{wrapfig}
        \usepackage{amsmath,amssymb,amsfonts}
        \usepackage{mathtools}
        \usepackage{wasysym}
        \usepackage{pgfplots}
        \usepackage{rotating}
        \usepackage{fancyvrb}
        \usepackage[short,nocomma]{optidef}
        \graphicspath{ {./images/} }

        \usepackage{algorithm,algpseudocode}
        \usepackage{lipsum}
        \usepackage{multirow}
        \usepackage{makecell}
        \usepackage{setspace}
        \usepackage{float}

\begin{document}
		
		\def\spacingset#1{\renewcommand{\baselinestretch}%
			{#1}\small\normalsize} \spacingset{1}
		
		\if0\blind
		{
			\title{Practice-Based Optimization for the\\
            Strategic Locomotive Assignment Problem}
			\author{Yunji Kim, Amira Hijazi, Kevin Dalmeijer, and Pascal Van~Hentenryck \\
			H. Milton Stewart School of Industrial and Systems Engineering, \\ Georgia Institute of Technology, Atlanta, GA, USA}
			\date{}
			\maketitle
		} \fi
		
		\if1\blind
		{

            \title{\bf \emph{IISE Transactions} \LaTeX \ Template}
			\author{Author information is purposely removed for double-blind review}
			
\bigskip
			\bigskip
			\bigskip
			\begin{center}
				{\LARGE\bf \emph{IISE Transactions} \LaTeX \ Template}
			\end{center}
			\medskip
		} \fi
		\bigskip
    
	\begin{abstract}
    \noindent This study addresses the challenge of efficiently assigning locomotives in large freight rail networks, where operational complexity and power imbalances make cost-effective planning difficult. It presents a strategic optimization framework for the Locomotive Assignment Problem (LAP), developed in collaboration with a major North American Class I Freight Railroad. The problem is formulated as a network-based integer program over a cyclic space-time network, producing a repeatable weekly locomotive assignment plan. The model captures a comprehensive set of real-world operational constraints and jointly optimizes the placement of pick-up and set-out locomotive work events, improving the effectiveness of downstream planning. To solve large-scale instances exactly for the first time, novel reduction rules are introduced to dramatically reduce the number of light travel arcs in the space-time network. Extensive computational experiments demonstrate the performance and trade-offs on real instances under a variety of practical constraints. Beyond delivering scalable, high-quality solutions, the proposed framework serves as a practical decision-support tool grounded in the operational realities of modern freight railroads. \\
	\end{abstract}
			
	\noindent%

	\spacingset{1.5} 


\section{Introduction}\label{sec:1}

Railroads play a vital role in intermodal freight transportation, offering unrivaled fuel-efficiency, cost-effectiveness, and safety among surface modes. As global freight volumes are projected to grow tremendously and policy initiatives promote a modal shift toward rail, improving the efficiency of rail operations has become more critical than ever \citep{forward302030, DOT_Bureau_of_Transportation_Statistics_2}. Central to this challenge are locomotives,  critical assets that require coordinated planning to ensure adequate power coverage across large, geographically dispersed networks.

Locomotive planning typically follows strategic, tactical, and operational levels  \citep{piu2015introducing}. Strategic decisions provide foundational input for subsequent planning phases considering shorter planning horizon and higher level of details. While the Operations Research (OR) community has made considerable progress on the Locomotive Assignment Problem (LAP) at the tactical and operational levels, strategic planning has received comparatively less attention. Yet, this level holds the greatest potential for long-term efficiency gains, with upstream decisions exerting strong influence over downstream decision quality.

This paper studies the \textbf{strategic Locomotive Assignment Problem} (LAP) that produces a repeatable, network-wide weekly plan that positions locomotives effectively across space and time. Specifically, it focuses on strategic network design for power planning and repositioning, motivated by a collaboration with a major North American Class I freight railroad. Unique challenges arise in topographically diverse regions, such as the mountainous Western U.S., where steep, variable grades cause substantial fluctuations in power needs over the course of a journey. These dynamics pose a greater impetus for intermediate power changes, referred to as \textit{work events} or \textit{power-changing points (PCP)} (e.g., \cite{ziarati1997locomotive}), which involve setting out or picking up locomotives mid-route. Unlike prior studies that assume work event locations as inputs, this study endogenizes these decisions, leading to enhanced system-wide outcomes.

In addition, disparities in train volume across the network lead to inherent imbalances in locomotive availability. While repositioning via scheduled trains (\textit{deadheading}) is effective, it is often inevitable to create new locomotive movements through \textit{light traveling}. Existing literature relies heavily on heuristic methods to generate a limited set of light travel arcs, such as those proposed in \cite{ahuja2005network}, without exploring the exact formulation using the full arc set. This study addresses that gap by incorporating the complete set of light travel arcs and applying reduction techniques to maintain computational tractability.

To ensure practical relevance, the proposed model incorporates a rich set of operational constraints and business rules, developed through extensive collaboration with industry stakeholders. In doing so, this work bridges the gap between academic optimization models and the realities of industrial deployment. The resulting framework offers not only cost-effective and strategically efficient plans, but also a practical decision-support tool grounded in the operational realities of a modern freight railroad.

\vspace{-0.5cm}
\subsection{Literature Review}
Locomotive planning has been a focal point in optimization studies for rail transportation since the mid-1970s. For comprehensive reviews on locomotive scheduling and related optimization models, readers are referred to \cite{cordeau1998survey}, \cite{ahuja2005network}, \cite{nemani2010or}, and \cite{piu2014locomotive}.
Early foundational works, such as \cite{florian1976engine}, introduced multicommodity network flow models, employing Bender's decomposition to tackle computational complexity. This approach was further extended in \cite{ziarati1999branch}, where the problem is decomposed into smaller time windows and solved by branch-and-cut integer programming algorithms \citep{ziarati1997locomotive}. The problem is modeled as a multicommodity flow problem on a space-time network. 

Since the minimum number of locomotives required to pull a train is more than one, many papers consider a heterogeneous fleet consisting of multiple locomotive types with different attributes, such as horsepower (HP). A popular aspect of LAP is to consider consists, formed of one or more locomotive types to meet HP requirements--for example, \cite{piu2015introducing} focused on preliminary consist selection. 
In the present study, locomotives are all of the same type and each train requires different integer numbers of them. This is an operationally reasonable assumption made by the industry partner at the strategic planning level, especially as they move towards a homogeneous fleet of high-horsepower.

Due to the complexity of freight rail systems and interdependencies among various decision-making components, decisions are made on three sequential levels as described in \cite{piu2015introducing}. Note that problem names and scopes may vary in related articles, not necessarily following the same classifications discussed here. 
At the strategic level, studies focus on high-level formulations that capture long-term goals like fleet sizing. \cite{powell2014locomotive} and \cite{bouzaiene2016single} develop a family of optimization models called the Princeton locomotive and shop management (PLASMA) system in collaboration with Norfolk Southern. A strategic variant PLASMA/SC is formulated as a deterministic single-commodity problem, where locomotives are of the same type. PLASMA/MC and PLASMA/MA versions correspond more closely to tactical and operational levels, tracking locomotives of heterogeneous types individually over the planning horizon. One important contribution was the use of approximate dynamic programming to address uncertainty.

The tactical problems delve into mid-term decision-making that aligns resources with short-to-mid-term operational objectives. For example, \cite{cordeau2000benders} employs Bender's decomposition to simultaneously assign cars and locomotives, optimizing coordination between these assets \citep{cordeau2001simultaneous}. 
\cite{ahuja2005solving} proposes integer linear programming (ILP) formulations based on space-time networks, and establishes the NP-Hard nature of LAP to underscore its complexity. 
More recent work such as \cite{scheffler2020mip} have considered the (dis-)connecting processes of locomotives and cars with heuristics to reduce the solution space. 
\cite{ortiz2021locomotive} builds upon the tactical LAP by incorporating the decision of distributed power in both locomotive-based formulations (LBF) and consist-based formulations (CBF). The computational performance was significantly improved using a Benders decomposition-based algorithm. 
\cite{frisch2021solving} solves the tactical level problem with maintenance constraints in collaboration with Rail Cargo Austria (RCA).

Lastly, the operational level addresses detailed, short-term planning tasks, including real-time adjustments and maintenance considerations. For instance, recent work of \cite{miranda2022time} focuses on time-sensitive locomotive scheduling, integrating repair shop routing into the planning process. Similarly, \cite{vaidyanathan2008real} examines real-time adjustments to operational plans, and \cite{vaidyanathan2008locomotive} extends this by incorporating maintenance constraints into the allocation process. These works differ in their emphasis on either immediate responsiveness or the integration of maintenance logistics.

\vspace{-0.25cm}
\subsection{Contributions}
\vspace{-0.1cm}
This paper advances strategic-level LAP by constructing a comprehensive space-time network that captures the weekly train schedule and reflects real-world operations. The main contributions are as follows.
\begin{itemize}
\item This paper presents an integer programming model that jointly optimizes locomotive assignments and the timing and location of work events. Unlike prior models with fixed work event locations, this framework optimizes them directly, improving input quality for downstream planning. The strategic LAP generates an assignment plan given the weekly train schedule and varying power requirements. 
\item This is the first work to solve the full LAP using the complete set of light travel arcs. To ensure tractability, reduction rules are introduced that exploit the temporal and spatial structure of \textit{cyclic} space-time networks, substantially reducing the number of arcs without loss of optimality. 
Despite the combinatorial complexity, the proposed approach enables global optimization at scale for large real-world instances. 
Notably, it is also the first to report true optimality gaps for heuristic methods in the literature. 
\item A wide range of operational constraints regarding locomotive work events are incorporated in order to enhance the practical applicability. These constraints, inspired by real use-cases, ensure that the proposed solutions are not only theoretically efficient but also realistic and deployable in the real world, bridging the gap between optimization and industrial practice.
\item Comprehensive computational experiments are conducted to assess the robustness of strategic LAP and its sensitivity to different cost parameters and operational constraints. These experiments offer valuable insights into how decision-makers can adjust the trade-offs between fleet size, repositioning flexibility, and operational overhead to achieve optimal outcomes under different scenarios.
\end{itemize}

The remainder of the paper is organized as follows. Section \ref{sec:2} provides a detailed description of the problem. Section \ref{sec:3} demonstrates how the problem is translated into a space-time network modeling framework, including its optimization formulation as well as methods for generating light travels. 
Sections \ref{sec:4} reports computational experiments and corresponding analyses. Conclusions follow in Section \ref{sec:5}.

\section{Problem Description}\label{sec:2}
The strategic LAP focuses on establishing a weekly plan for assigning locomotives to scheduled trains. Since the train schedule is assumed to be repeated every 7-day period, the locomotive assignment plan must be repeatable every week as well. In the following section, the main inputs, operational constraints as well as the costs considered in the problem are described. The main output involves the assignment of locomotives throughout the week as well as the scheduling of locomotive work events such as pick-ups and set-outs. The decisions made at the strategic level, including the fleet size and work event locations, form the foundation for the decision-making process in tactical and operational problems.

\subsection{Problem Data} \label{sec:2.1}
The primary inputs to the LAP include a rail network and a weekly train schedule that specifies the arrival and departure times, inbound and outbound terminals, and power demands of each train. A train consists of a group of railcars and locomotives. The number of locomotives required for a train is determined by the total weight of the railcars, as well as the operating conditions along the route. 
Importantly, railcars can be attached or detached at intermediate terminals during a train’s journey. This implies the train’s weight and composition may vary along its route, which in turn affects the motive power needed at different segments. As a result, the demand for locomotives is not static, but shaped by multiple operational factors including load, terrain, and dynamic changes in train makeup.

This section outlines the key components of the input data and highlights how they interact to influence locomotive assignments. 

\vspace{-0.3cm}
\paragraph*{Rail Network} The rail network is represented as a set of terminals $K$ connected by rail tracks. The connectivity between terminals is characterized by transit times, which may be asymmetric due to varying ruling grades and terrain. Transit time between terminals $i,j \in K$ is denoted as $\delta(i,j)$.

\vspace{-0.3cm}
\paragraph*{Train Schedule} 
A weekly train schedule outlines movements of trains $T$ across the network, detailing origins, destinations, and intermediate stops. Each train $t \in T$ has its journey divided into $s_t$ train legs, indexed by $I_t = \{1, \dotsc, s_t\}$, representing sequential segments. The intermediate stops occur from operations such as railcar work events (pick-up/set-out), crew changes, fueling, or inspections. Each leg is defined by its origin and destination terminals and scheduled arrival and departure times. 
The complete set of train legs is denoted by $A_T$, with each leg identified by its train $t \in T$ and sequence index $i \in I_t$.

\vspace{-0.3cm}
\paragraph*{Locomotives} The fleet of locomotives is assumed to be homogeneous, meaning all units have identical pulling capability.
This assumption aligns with current practice at the railroad industry partner for strategic planning.
For each train leg $l \in A_T$, the required motive power is pre-calculated as $b_l$, representing the integer number of locomotives needed to sufficiently pull the train's tonnage over the leg's ruling grade. The value $b_l$ reflects
both the terrain and dynamic changes in train weight due to railcar work events.

\vspace{-0.3cm}
\paragraph*{Railcar Work Events} During a train's journey, railcars may be attached (pick-up) or detached (set-out) at intermediate terminals, altering the train's tonnage and thus the motive power requirements. In some cases, dropping off railcars and picking up new ones can occur simultaneously. 
To capture this, binary parameters $r^{pu}, r^{so}, r^{no}$, and $r^{both}$ are defined to indicate whether a pick-up, set-out, neither, or both types of railcar work events occur. These values are defined for every intermediate stop. 

These events are particularly relevant as they coincide with opportunities to adjust the number of locomotives for upcoming train legs--referred to as \textit{locomotive work events}. The alignment of railcar and locomotive work events impacts locomotive work events costs.
Unlike current plans in practice that use a limited set of locations, this paper considers all intermediate stops as potential sites for locomotive work events.

\subsection{Flow Balance and Locomotive Repositioning} \label{sec:2.2}
In this study, locomotives are classified as \emph{active} (engine on), or \emph{inactive} (engine off). Active locomotives correspond to the input $b_l$ for each train leg $l \in A_T$. To ensure active locomotives are available where and when needed, locomotives are repositioned across the network while inactive. Since locomotives are homogeneous and not tracked individually, what matters is having the right number of units in the right place at the right time—not which specific ones.

At each terminal, the weekly locomotive balance is determined by horsepower (HP) inflow from inbound trains versus outflow for outbound trains. A terminal is \emph{balanced} if its locomotive inventory at the start of the week matches that at the end. However, this equilibrium is rare due to inherent inequalities in supply and demand. Terminals with greater supply experience a \emph{surplus}, resulting in excess idle locomotives, while those with greater demand face a \emph{deficit}, causing shortages for outbound trains. These imbalances necessitate strategic repositioning of inactive locomotives for a balanced plan. 

Locomotive repositioning occurs via two methods: \emph{deadheading} (DH) or \emph{light traveling}. In DH, inactive locomotives are attached to scheduled trains, leveraging existing movements which makes this a cost-efficient option. This approach naturally creates locomotive work events. The locations and number of locomotives to be picked up or set out at each stop are critical in strategic LAP. Pick-ups require preparing locomotives at terminals prior to departure, while set-outs involve detaching units that must undergo mandatory maintenance and inspection. Afterward, they enter ground inventory as idle until reassigned. The durations of both activities are user-defined. 
Alternatively, a group of inactive locomotives can form a \emph{light travel} train, with only the lead locomotive pulling. Although substantially more expensive due to dedicated crew for empty moves lacking productive freight movement, it is often necessary to maintain network balance.

\subsection{Costs} \label{sec:2.3}
The total cost of a weekly LAP includes locomotive ownership, work events, light travel, and relocation. 
Locomotive work events, such as picking up or setting out locomotives, often occur at the same intermediate stops where railcar work events take place. These two activities can either align or conflict in their type, which significantly affects the cost of locomotive work events. The lowest cost, $c_1$, applies when locomotive and railcar work events are aligned—e.g., both pick-ups or both set-outs—allowing joint execution with minimal disruption. If the events are mismatched—e.g., a locomotive pick-up alongside a railcar set-out—a higher cost $c_2$ is incurred due to the need for sequential handling. The highest cost $c_3$ arises when a locomotive event occurs at a stop with no railcar event at all, requiring a stand-alone activity with the greatest operational burden. Costs are incurred per locomotive event, regardless of the number of units involved. A special case occurs when both railcar pick-ups and set-outs are scheduled at a terminal; in this case, either type of locomotive event is treated as aligned and costs $c_1$. This cost structure incentivizes aligning locomotive repositioning with existing railcar activity to reduce overhead and improve efficiency.

For light travel, the cost consists of a fixed component $e_l$ and a variable component $g_l$ for each light travel arc $l \in A_L$.
Relocation cost $g_l$ depends on the number of inactive locomotive units moved and length (transit time) of $l$, and is applied to both light travel and DH.
The fixed charge $e_l$ reflects the cost for assigning additional crew, and it is based on the length of $l$.
All light travel arcs between the same terminal pair share identical costs regardless of departure time.

Finally, a weekly ownership cost $q$ is applied for each locomotive unit. 
The optimization model aims to minimize the sum of these costs while adhering to a variety of operational constraints.
For reference, a summary of the cost parameters is included in Table~\ref{table:Nomenclature}.

\vspace{-0.3cm}
\section{Modeling Framework} \label{sec:3}
\vspace{-0.3cm}
The strategic LAP is modeled as a single-commodity flow problem on a space-time network that represents both the physical locations and the time at the minute level.
Figure~\ref{fig:tikzfigure_network} will serve as a running example throughout this section.
Let $G=(N,A)$ be a graph with nodes $N$ and arcs $A$. Each arc represents a locomotive activity from one space-time node to another, with the flow on an arc denoting the number of locomotives involved in that movement. Table \ref{table:Nomenclature} summarizes the notations used throughout this paper. The following sections detail the space-time network and formally introduce the optimization model to solve the strategic LAP.

\begin{table}[!]
    \centering
    \caption{Nomenclature used in this Study}
    \begin{singlespace} 
    \resizebox{0.9\linewidth}{!}{
    
    \begin{tabular}{l l }%
    \hline
    \textbf{Symbol} & \textbf{Definition} \\ 
    \midrule    

    \textbf{Sets}: &  \\
        $T$& Set of scheduled trains \\
        $K$& Set of terminals \\
        $N, A$& Set of nodes and arcs \\
        $N_D, N_A$& Set of departure and arrival nodes \\
        $N_G = N_I \cup N_E \cup N_R $& Set of ground nodes (initial, ground-departure, and arrival-ground)\\
        $A_T$& Set of train arcs\\
        $A_C$& Set of train-leg transition arcs\\
        $A_E, A_R$& Set of ground-departure and arrival-ground arcs\\
        $A_{PU} \subset A_E, A_{SO} \subset A_R$&Set of set-out and pick-up arcs\\
        $A_L$& Set of light travel arcs\\
        $I[n], O[n]$& Sets of inbound and outbound arcs for each space-time node $n \in N$ \\
        $S$& Set of arcs wrapping around the time horizon \\ 
    
    \textbf{Parameters}: & \\
        $\delta(i,j)$& Transit time between terminals $i,j \in K$ in minutes \\
        $s_t$ & Number of train legs in train $t \in T$\\
        $I_t$ & Set of sequence indices $[1, s_t]$ for legs in train $t \in T$ \\
        $b_{l}$ & Number of active locomotives for train leg $l \in A_T$ \\
        $f$& Maximum number of total locomotives per train leg\\
        $\rho^U$& Maximum number of locomotives per light travel train\\
        $r_{l}^{so}, r_{l}^{pu}, r_{l}^{no}, r_{l}^{both}$& Binary parameter indicating railcar work events occurrence at transition arc $l \in A_C$\\ 
        $c_1, c_2, c_3$& Per-event cost of locomotive work event varied by type\\
        $q$& Per-unit weekly locomotive ownership cost\\ 
        $e_l$& Fixed cost of activating a light travel arc $l \in A_{L}$\\
        $g_l$& Per-unit cost of assigning an inactive locomotive on arc $l \in A_{T} \cup A_{L}$ \\
    
    \textbf{Decision Variables}: & \\
    $x_{l}$ & Integer variable indicating the unit of locomotives on $l \in A$\\
    $y_{l}^{so}, y_{l}^{pu}$ & Binary variables taking value 1 if arc $l \in A_{SO} / A_{PU}$ is used and zero otherwise, resp.\\
    $u_{l}$ & Integer variable indicating the number of light travel trains on arc $l \in A_{L}$\\   
    
    \bottomrule
    \end{tabular} 
    }
    \end{singlespace}
    \label{table:Nomenclature}
    \vspace{-0.4cm}
\end{table}

\vspace{-0.3cm}
\subsection{Space-Time Network} \label{sec:3.1}
The space-time network models locomotive activities through nodes and arcs. Each node $n \in N$ is defined by its spatial (terminal) and temporal attributes, as well as its type. The node set $N$ includes train departure nodes $N_D$, arrival nodes $N_A$, and ground nodes $N_G$. Locomotives on ground are stationary at a terminal as inventory, preparing to be attached to a train, or having just been removed from service. Accordingly, ground nodes include three types. \textit{Initial} node exists at each terminal at time zero which form a set $N_I$. \textit{Ground-departure} nodes $N_E$ indicate transitions from ground to train service, and \textit{arrival-ground} nodes $N_R$ mark when set-out locomotives join ground inventory after maintenance.

Arcs $A$ represent locomotive activities such as train travel, work events, or ground idle time.
\textit{Train} arcs $A_T$ correspond to scheduled train legs, connecting departure and arrival nodes. Each arc $l \in A_T$ is associated with a train $t \in T$ and a leg $i \in I_t$. These are the only arcs on which the active locomotives flow.

\begin{figure}[!]
  \centering
  \resizebox{0.8\textwidth}{!}{\begin{tikzpicture}[
    node/.style={circle, minimum size=6mm, inner sep=0pt},
    graynode/.style={node, fill=gray!70},
    rednode/.style={node, fill=red!40}, 
    bluenode/.style={node, fill=blue!30}, 
    purplenode/.style={node, fill=purple!80}, 
    yellownode/.style={node, fill=orange!30}, 
    thickarrow/.style={->, thick},
    verythickarrow/.style={->, ultra thick},
    dashedarrow/.style={->, very thick, dashed},
    dottedarrow/.style={->, dotted, very thick},
]

\def\yA{0}     
\def\yB{-4}    
\def\yC{-8}    

\foreach \y/\label in {\yA/Terminal 1, \yB/Terminal 2, \yC/Terminal 3} {
    \draw[thick] (-1,\y) -- (12,\y);
    \node[anchor=east] at (-1.2,\y) {\textbf{\label}};
}

\node[graynode] (g1) at (-0.5,\yA) {};
\node[purplenode] (p1) at (11.5,\yA) {};
\node[below=2pt of p1] {\scriptsize arvl-grd node};

\draw[thickarrow] (g1) -- (p1);
\draw[dottedarrow, bend right=20] (p1) to (g1);

\node[graynode] (g2) at (-0.5,\yB) {};
\node[below=2pt of g2] {\scriptsize initial node};

\node[yellownode] (y1) at (4,\yB) {};
\node[purplenode] (p2) at (8.5,\yB) {};

\node[rednode] (r1) at (6,-3) {}; 
\node[bluenode] (b1) at (9,-1) {}; 
\node[below=2pt of b1] {\scriptsize arvl node};

\draw[thickarrow] (g2) -- (y1);
\draw[dashedarrow] (y1) -- (r1) node[near end, left, xshift=-5pt] {\scriptsize pick-up arc};
\draw[verythickarrow] (r1) -- (b1);
\draw[thickarrow] (b1) -- (p1);

\draw[thickarrow] (y1) -- (p2);
\draw[dottedarrow, bend right=40] (p2) to (g2);
\node[graynode] (g3) at (-0.5,\yC) {};
\node[yellownode] (y2) at (1,\yC) {};
\node[below=2pt of y2] {\scriptsize grd-dep node};

\node[rednode] (r2) at (2,-7) {}; 
\node[right=2pt of r2] {\scriptsize dep node};

\node[bluenode] (b2) at (4.5, -5) {}; 

\draw[thickarrow] (g3) -- (y2);
\draw[thickarrow] (y2) -- (r2);
\draw[verythickarrow] (r2) -- (b2) node[near end, left, xshift=-5pt] {\scriptsize train arc};
\draw[thickarrow] (b2) -- (r1) node[pos=0.7, right] {\scriptsize transition arc};

\draw[dashedarrow] (b2) -- (p2) node[midway, below, yshift=-5pt] {\scriptsize set-out arc};

\draw[dottedarrow, bend left=40] (y2) to (g3);

\end{tikzpicture}}
  \caption{An example space-time network}
  \label{fig:tikzfigure_network}
  \vspace{-0.3cm}
\end{figure}
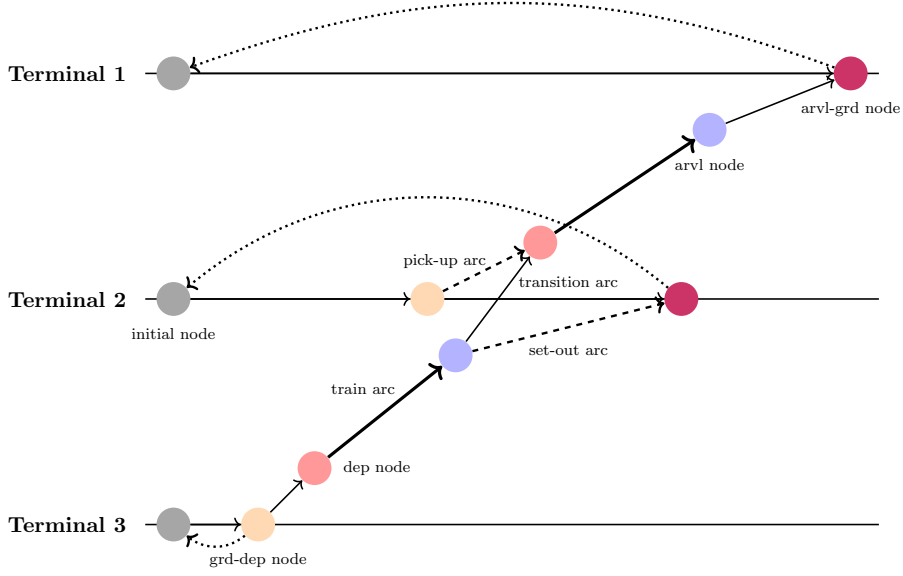

Before traveling on train legs, locomotives transition from ground via \textit{ground-departure} arcs $A_E$. They connect ground-departure nodes $N_E$ to departure nodes $N_D$, and model train preparation and building. 
A subset of $A_E$, excluding arcs associated with first train legs ($\{1 \in I_t: \forall t \in T\}$), forms the set of \textit{pick-up} arcs $A_{PU}$. These arcs represent pick-up decisions at intermediate stops and exclude initial departures which must always be used.

Upon arriving at an intermediate stop after completing a train leg, locomotives may either remain attached for the next leg or be set out. If continuing, they flow through \textit{train-leg transition} arcs $A_C$, which connect the arrival node of $i$-th train leg to the departure node of the $(i+1)$-th leg. 
These arcs model dwell time at the shared stop between consecutive legs, and ensure continuity of the journey. Such arcs exist at all intermediate stops for trains with multiple legs, and the aforementioned $r$ railcar parameters are defined on them.
Alternatively, if locomotives are set out, they flow along \textit{arrival-ground} arcs $A_R$, linking arrival nodes $N_A$ to corresponding arrival-ground nodes $N_R$. These arcs represent the set-out operations and mandatory inspection period locomotives undergo before becoming eligible for train assignment. The subset of $A_R$ excluding the arcs following the final train legs forms the \textit{set-out} arcs $A_{SO}$, which designate where set-out decisions are made. In addition, \textit{light travel} arcs $A_L$ between ground nodes $N_G = N_I \cup N_E \cup N_R $ are pre-generated and added to the space-time network. Section~\ref{sec:3.3} discusses this process in detail.

At each terminal, \textit{ground} arcs $A_G$ connect ground nodes sorted by time sequentially and represent idle locomotives. 
To model a cyclic weekly horizon, a subset $S \subseteq A_G$ defines \emph{wrap-around} arcs, connecting the final ground node at each terminal back to its initial node. They capture locomotive flow transitions from the end of one week to the start of the next, ensuring inventory continuity across weeks.
Furthermore, arcs whose time intervals cross the horizon (e.g., Friday to Monday of the following week) have their time attributes adjusted modulo the horizon length, to preserve cyclic structure. Collectively, arcs in $S$ model a time-wrapped network, a standard technique in scheduling (e.g., \cite{ahuja2005network}).
Finally, for each node $n \in N$, the sets $I[n]$ and $O[n]$ define its inbound and outbound arcs, respectively.

Figure \ref{fig:tikzfigure_network} serves as an example space-time network of a train journey with two train arcs represented with bold arcs. The dashed arcs represent the pick-up and set-out arcs, and the dotted arcs represent the cyclic behavior. Upon arriving at terminal 2, locomotives arriving on the first train arc may either be dropped off using the set-out arc or continue onto the next leg using the transition arc. Note that additional locomotives will be picked up only if none are set out. That is, at most one of the pick-up or set-out arcs can be used to change the total number of locomotives between train legs.

\vspace{-0.4cm}
\subsection{Formulation}\label{sec:3.2}
\vspace{-0.2cm}
The strategic LAP is modeled through integer flows on the space-time network. Each arc $l \in A$ is associated with an integer decision variable $x_l$, representing the number of locomotives assigned to it.
For each pick-up arc $l \in A_{PU}$, a binary variable $y_l^{pu}$ is defined, taking value 1 if $x_l > 0$ (i.e., the arc is used), and 0 otherwise. Similarly, for each set-out arc $l \in A_{SO}$, a binary variable $y_l^{so}$ is defined. These fixed charge variables model the occurrence of locomotive work events, while the associated $x_l$ variables capture the actual number of locomotives picked up or set out. 
Lastly, let $u_l$ be an integer variable indicating the number of light travel trains operating on arc $l \in A_L$.

The space-time network and variable definitions directly lead to the compact Model~\eqref{form:base}.
Objective \eqref{obj_1} minimizes the total cost consisting of four terms.
The first term determines the total weekly locomotive ownership cost.
Since the network is time-wrapped, flows on the set of wrap-around arcs $l \in S$ capture locomotives transitioning between consecutive weeks, and count the fleet size over the planning horizon (e.g., \cite{ahuja2005network, ortiz2021locomotive}).
The second term represents the total cost of DH, capturing the expense of attaching inactive locomotives ($x_l$ - $b_l$) to scheduled trains.
The third term captures costs of activating light travel arcs and repositioning locomotives through them.
Finally, the function \( RC(y) \) penalizes locomotive work events based on railcar work event parameters \( r_l \) and the associated costs (Section~\ref{sec:2.3}). For each transition arc \( l = (i, j) \in A_C \), the penalty depends on whether a locomotive set-out occurs at node \( i \in N_A \) or a pick-up at node \( j \in N_D \), and this is captured by  binary variables \( y_{l'}^{so} \) and \( y_{l^*}^{pu} \), respectively. Each \( r_l \) is weighted by cost coefficients \( c_1, c_2, c_3 \) to compute the total penalty:
\begin{multline*}
RC(y) = \sum_{l = (i,j) \in A_C} \, \sum_{l^* \in A_{PU} \cap I[j]} \, \sum_{l' \in A_{SO} \cap O[i]} \Big[ 
r_l^{so} (c_1 y_{l'}^{so} + c_2 y_{l^*}^{pu}) 
+ r_l^{pu} (c_2 y_{l'}^{so} + c_1 y_{l^*}^{pu}) \\
+ r_l^{no} (c_3 y_{l'}^{so} + c_3 y_{l^*}^{pu}) 
+ r_l^{both} (c_1 y_{l'}^{so} + c_1 y_{l^*}^{pu}) 
\Big].
\end{multline*}

\begin{figure}[t]
\vspace{-0.2cm}

\begin{mini!}
%
	{}
%
	{\sum_{l \in S} q x_l + \sum_{l \in A_T} g_l (x_l - b_l) + \sum_{l \in A_L} (g_l x_l + e_l u_l) + RC(y) \label{obj_1}}
%
	{\label{form:base}}
%
	{}
%
%
	\addConstraint
	{b_l \le x_l}
	{\le f}
	{\forall l \in A_T \label{cstr:cap}}
	\addConstraint
	{\sum_{l \in I[n]} x_l}
	{= \sum_{l \in O[n]} x_{l} \qquad}
	{\forall n \in N \label{cstr:flow_1}}
	\addConstraint
	{x_l}
	{\le f y_l^{so}}
	{\forall l \in A_{SO} \label{cstr:so_1}}
	\addConstraint
	{x_l}
	{\le f y_l^{pu}}
	{\forall l \in A_{PU} \label{cstr:pu_1}}
	\addConstraint
	{x_l}
	{\le \rho^U u_l}
	{\forall l \in A_L \label{cstr:emoveub_1}}
 	\addConstraint
	{x_l}
	{\in \mathbb{N}_0}
	{\forall l \in A \label{cstr:var1_1}}
 	\addConstraint
	{y_{l}^{so}}
	{\in \{0, 1\}}
	{\forall l \in A_{SO} \label{cstr:var2_1}}
 	\addConstraint
	{y_{l}^{pu}}
	{\in \{0, 1\}}
	{\forall l \in A_{PU} \label{cstr:var3_1}}
 	\addConstraint
	{u_l}
	{\in \mathbb{N}_0}
	{\forall l \in A_L \label{cstr:var4_1}}
\end{mini!} 

\caption{Optimization Model for the Strategic LAP}
\vspace{-0.2cm}
\end{figure}

Constraints~\eqref{cstr:cap} enforce that the flow on each train arc satisfies the minimum power requirement, yet does not exceed the maximum number of $f$ locomotives per train.
Constraints \eqref{cstr:flow_1} impose flow conservation on every space-time node, enforcing a cyclic solution that may be repeated every week. Constraints \eqref{cstr:so_1} and \eqref{cstr:pu_1} activate set-out and pick-up arcs when they are used by locomotives. 
Constraints~\eqref{cstr:emoveub_1} ensure that light travel arcs can only be used when light trains are operated.
Furthermore, they limit the number of locomotives per light travel train to $\rho^U$.
Finally, the variables are defined by Equations~\eqref{cstr:var1_1}-\eqref{cstr:var4_1}.
Note that the $u$-variables are integer instead of binary, such that multiple light trains may be sent over the same arc (at a cost).

\subsection{Generating Light Travel Arcs} \label{sec:3.3}
Light travel is often necessary to achieve power balance across the network, as relocating locomotives solely through deadheading is neither efficient nor sufficient. In principle, locomotives can light travel between any pair of terminals at each discretized time unit. However, this generates a vast number of candidate arcs, dramatically increasing problem complexity. To ensure tractability, prior studies have relied on heuristics to select a manageable set of light travel arcs $A_L$. 
This section presents a method to reduce the full set of candidate light travel arcs and solve the exact problem without loss of optimality. It also reviews a widely used Minimum Cost Flow (MCF) approach from the literature. 

\vspace{-0.3cm}
\subsubsection{Exact Method} \label{sec:3.3.1}
\vspace{-0.2cm}
The complete set $\bar{A_L}$ contains all possible light travel arcs within the space-time network over the planning horizon. While optimizing over this complete set yields an exact solution, the resulting model is computationally intractable due to the sheer number of variables and constraints introduced.
To address this, a reduced subset $A_L \subseteq \bar{A_L} $ is constructed that preserves optimality while drastically shrinking the network size.

The key idea is to generate arcs only between existing ground nodes tied to discrete events (train arrivals and departures), rather than between arbitrary time points. This step effectively trims $\bar{A}_L$ without introducing new nodes, and avoids modeling unnecessary possibilities that would never appear in an optimal solution. Specifically, light travel arcs are created from arrival-ground nodes $N_R$ (where locomotives become idle) to ground-departure nodes $N_E$ (where they are required). 
This ensures that light travel occurs only within meaningful transfer windows—beginning immediately when locomotives become available and ending just before they are reassigned. The reduction proceeds in two steps:
\begin{enumerate}
    \item Earliest Reachability: Each arrival-ground node is connected to the first available ground-departure node at each destination terminal that it can reach. Since locomotives can wait idly on the ground at no cost once they arrive, reaching the destination as early as possible suffices.

    \item Latest Origin Filtering: If multiple arrival-ground nodes from the same origin terminal connect to the same ground-departure node, only the arc from the latest such origin node is retained. Earlier departures to the same destination would result in idle time at the destination and thus offer no advantage.
\end{enumerate}

\noindent The resulting arc set $A_L$ may be significantly smaller and does not sacrifice optimality.
This can be seen from the fact that the fixed charge and variable cost for light travel are constant for arcs between the same terminals.
Hence, consolidating light travel on fewer arcs can never increase the objective value.
Applying the reduction is straightforward, but some care must be taken in handling wrap-around arcs.
Details are provided in the Appendix.

\subsubsection{Minimum Cost Flow Heuristic} \label{3.3.2}
A common heuristic in locomotive planning is to generate light travel arcs by solving a minimum-cost flow (MCF) problem on a space network \citep{ahuja2005solving}. Terminals are classified as power sources or sinks based on the weekly net inflow and outflow of required locomotives: sources have excess outgoing flow, whereas sinks have excess incoming demand. The MCF is solved to determine optimal flows between each source-sink pair, with the cost coefficient $e_{ij}$ for light travel from terminal $i$ to $j$ defined as:

\begingroup
\singlespacing
\[
\renewcommand{\arraystretch}{0.9}
e_{ij} = \begin{cases}
    \delta_{ij} &\qquad \text{if } o_{ij} \leq 2 \\
    \delta_{ij}  \cdot \alpha &\qquad 2 < o_{ij} < \alpha\\
    \delta_{ij} \cdot \alpha^2 &\qquad \text{otherwise,}
\end{cases}
\]
\endgroup

\noindent where $\delta_{ij}$ is the railroad distance, $o_{ij}$ is the total number of trains operating between $i$ and $j$ in the weekly schedule, and $\alpha$ is a penalization factor discouraging flow through heavily serviced terminal pairs. If the optimal flow on an arc from solving MCF exceeds a pre-defined threshold, corresponding light-travel arcs are added to the space-time network. 
To include the time dimension, the planning horizon is partitioned into equal-length time windows. One space-time arc is inserted for each window, originating from an origin ground node and entering the first available ground node at the destination after travel time. If no eligible node exists within a window, the nearest neighboring window is used. 

While this heuristic effectively generates a subset of light-travel arcs to make the problem feasible, it provides no guarantee of selecting the optimal light-travel arcs. Moreover, adding arcs at uniform temporal intervals may lead to inefficiencies both in terms of computational effort and solution quality. 

\vspace{-0.1cm}
\subsection{Practice-based Model Extension: Work Event Restrictions} \label{sec:3.4}
\vspace{-0.1cm}
The optimization model in Section \ref{sec:3.2} assumes unconstrained work events scheduling, i.e., any terminal, any day, with no capacity limits. However, in the real-world, work event feasibility hinges on crew constraints, track and terminal capacity, and coordination among various stakeholders including field and dispatch teams.
To bridge this gap, this section introduces model extensions that incorporate practice-based constraints. These extensions explore constrained flexibility within real-world limits to produce a plan that is practical and executable.

\begin{table}[t]
    \centering%
    \caption{Additional Sets and Parameters used in Model Extensions}
    \resizebox{0.9\linewidth}{!} {
    \begin{tabular}{l l }%
    \hline
    \textbf{Symbol} & \textbf{Definition} \\ 
    \midrule   
    $D$& Set of days in the planning horizon (e.g., days of the week) \\
    $K^I$& Set of terminals without work events under the current plan\\
    $KD$& Set of feasible terminal-day pairs \\
    ${KD}^I$& Set of terminal-day pairs without work events under the current plan \\
    $\alpha^C, \alpha^D, \alpha^E, \alpha^F$& Maximum number of additional work event activations allowed per model version \\
    $\theta$& Upper limit on the number of work events allowed at a terminal on any given day \\
    $h_{k,d}$& Number of work events at terminal $k$ on day $d$ under the current plan \\
    $\lambda$& Permitted increase in work events at already active terminal-day pairs \\
    $\mathcal{L}_{kd}$& Set of train connection arcs $A_C$ taking place at terminal $k$ on day $d$ \\
    \bottomrule
    \end{tabular} }
    \label{table:addParams}
\end{table} 

The current locomotive work event schedule provided by the industry partner serves as a baseline. Table \ref{table:addParams} summarizes the additional notations used in the model extension.
Let $D$ be the set of discretized days in the planning horizon, and $KD$ the set of all feasible terminal-day combinations. The parameter $h_{k,d}$ gives the number of work events scheduled at terminal $k$ on day $d$ according to the current plan, and $\theta$ specifies the maximum number of allowable work events at a terminal on any given day. 

A terminal-day is \textit{active} if it has at least one work event; a terminal is \textit{active} if any of its days are active during the planning horizon. 
$K^I \subseteq K$ is the set of baseline inactive terminals with no scheduled work events across all days, i.e., $h_{k,d} = 0$ for all $d \in D$. Similarly, ${KD}^I \subseteq {KD}$ is the set of baseline inactive terminal-day pairs $(k,d)$ where $h_{k,d} = 0$. 

Five model extension versions are introduced and they are categorized into three groups by their level of deviation from the current plan: 
\begin{enumerate}
    \item \emph{Capacity Increase at Baseline-Active Terminal-Days (V1): }  \emph{V1} increases work event capacity at base-line active terminal-day pairs (${KD} \setminus {KD}^I$) by $\lambda$ units. No new activations are introduced.
    \item \emph{Incremental Expansions (V2, V3) : } These models extend the baseline gradually by introducing limited new activities. \emph{V2} allows up to $\alpha^C$ new terminals (from $K^I$) to be activated. If a terminal is active, it can have work events any day of the week.
    \emph{V3} allows up to $\alpha^D$ new terminal-day pairs (from ${KD}^I$), offering temporal control.
    \item \emph{Network Redesigns (V4, V5):} These versions disregard the baseline and optimize from scratch.
    \emph{V4} selects a new set of active terminals (up to $\alpha^E$) that can have work events any day.
    \emph{V5} activates any terminal-day pairs up to $\alpha^F$, providing the most flexible configuration. 
\end{enumerate}

These frameworks allow decision-makers to balance operational gains and implementation feasibility. Control parameters $\alpha$, $\lambda$, and $\theta$ regulate the scale and structure of changes. Detailed model formulations are provided in the Appendix.

\vspace{-0.2cm}
\section{Computational Experiments} \label{sec:4}
\vspace{-0.2cm}
This section presents computational results on real data to assess the empirical performances of strategic LAP formulation with practice-based constraints and the reduced, full set of light travel arcs. The optimization models were implemented in Python 3.9 using FICO Xpress 9.5.0, with four threads allocated to each instance. A termination criterion of an 8-hour time limit was applied unless specified otherwise. All experiments were conducted on a server equipped with dual-socket Intel Xeon Gold 6226 CPUs, each featuring 24 cores running at 2.7GHz.

\subsection{Benchmark Instances} The computational experiments were performed using real-world instances from historical data provided by the industry partner. Due to the complexity and scale of network operations, as well as the fragmented nature of the raw data, significant effort was devoted to data cleaning and preprocessing.
This process involved extensive validation to ensure the resulting data are accurate and suitable for the scope of this study. 
As part of the preprocessing, individual stations are aggregated into terminal-level representations based on geographic proximity, reflecting the operational reality that locomotives can be freely shared within a terminal. Furthermore, the data is restricted to particular trains that operate with high-horsepower (HHP) locomotives, which is the primary fleet of our case study.  
A typical weekly train schedule has over 10,000 train legs in the rail network comprising of hundreds of terminals. The resulting space-time network typically contains over 40,000 nodes and 60,000 arcs, and the corresponding optimization model has on the order of 90,000 decision variables and constraints, without counting the light travel arcs.

\subsection{Computational Performance of Exact vs. MCF} When implementing the exact method described in Section~\ref{sec:3.3}, over a million light travel arcs are generated: this is a dramatic reduction from the original candidate set, yet significantly larger than the hundreds of arcs produced by the MCF approach. As a result, the exact method yields millions of decision variables and constraints, making the problem size substantially larger. The MCF is solved with $\alpha$ as the mean of $o_{ij}$ (number of trains scheduled between terminals $i$ and $j$), and light travel arcs are added to the space-time network in every 8-hour window for OD pairs with optimal flow greater than 1. 

\begin{figure}[!]
    \centering
    \subfigure[True Optimality Gaps Over Time]{%
        \includegraphics[width=0.48\textwidth]{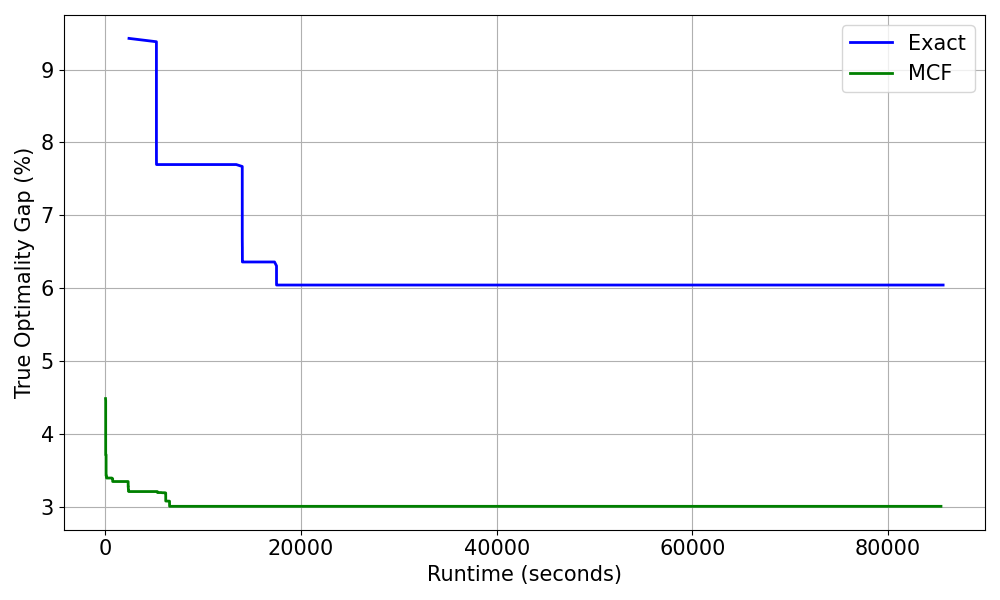}
        \label{fig:upper_lower_bounds}
    }
    \hfill
    \subfigure[Upper and Lower Bounds Over Time]{%
        \includegraphics[width=0.48\textwidth]{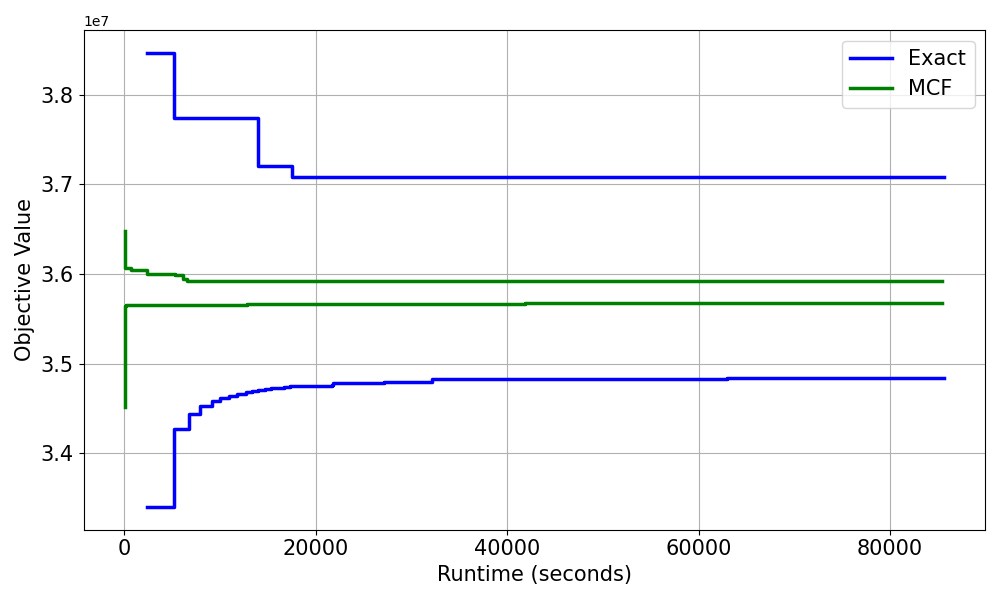}
        \label{fig:obj_gap_plot}
    }
    \caption{Performance Comparison between Light Travel Generation Methods}
    \label{fig:exact_MCF_comparison}
\end{figure}

Figure~\ref{fig:exact_MCF_comparison} compares the computational performance of the two methods with an 8-hour solver time limit. All optimality gaps are calculated with respect to the true lower bound from the exact method. The MCF-based model reaches a true gap of around $3\%$ within two hours and shows limited improvement thereafter. 
In contrast, the exact model requires a couple thousand seconds for presolve and does not reduce the gap below $6\%$ within the allotted time. The typical gap difference between the methods ranges from 2 to 4 percentage point. Despite the computational burden, the exact method offers key theoretical advantages: it guarantees convergence to the true optimum by solving the exact problem. It is also valuable for validating solution quality through a valid lower bound.
Most notably, \emph{calculating the true optimality gap of MCF in Figure~\ref{fig:upper_lower_bounds} would not be possible without the exact method}.

\subsection{Trade-Off Analysis: System Behavior under Different Cost Configurations}
In practice, decision makers often face conflicting objectives. For instance, reducing fleet size lowers ownership costs but requires more repositioning and work events. Improving one metric may therefore come at the expense of others. To better understand these trade-offs, a sensitivity analysis is conducted to scrutinize how varying cost parameters influence assignment plan characteristics. 

\begin{figure}[!]
    \centering
    \subfigure[Varying Ownership Cost ($q$) \label{fig:sens_fleet}] {\includegraphics[width=\textwidth]{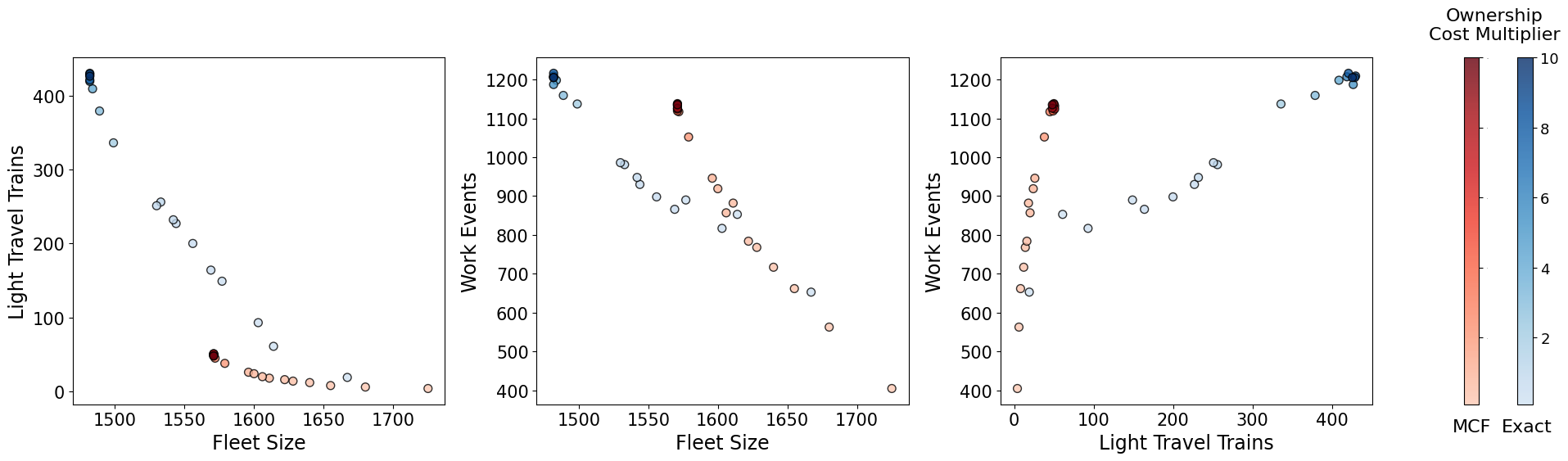}} \\
    \subfigure[Varying Light Travel Crew Cost ($e$) \label{fig:sens_lt}]{\includegraphics[width=\textwidth]{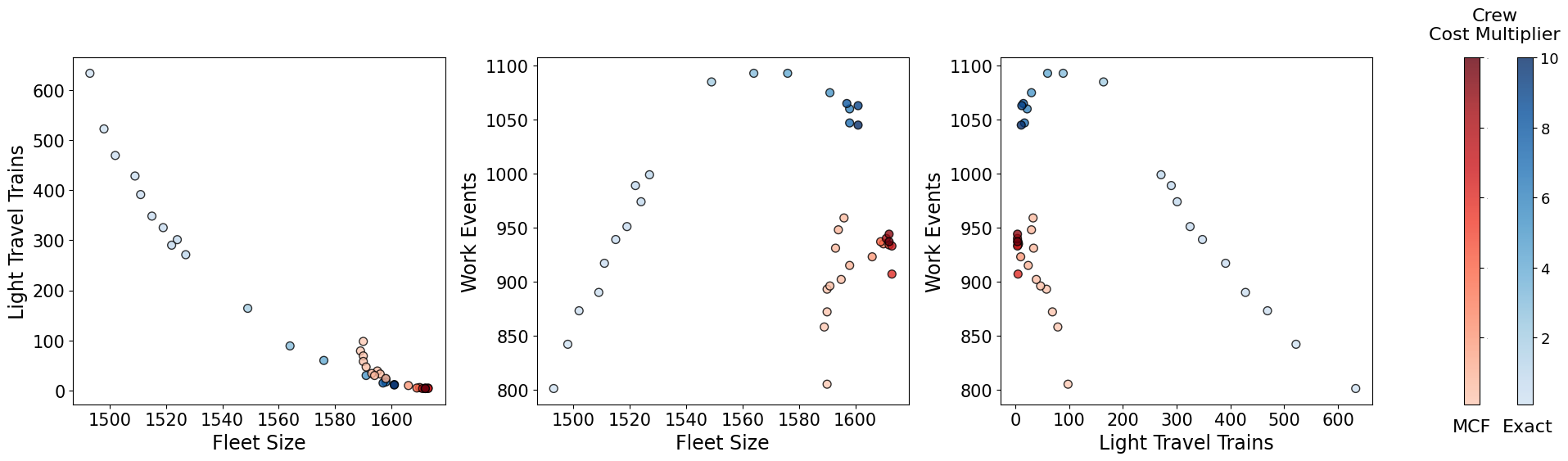}} \\
    \subfigure[Varying Work Events Costs ($c$) \label{fig:sens_wc}]{\includegraphics[width=\textwidth]{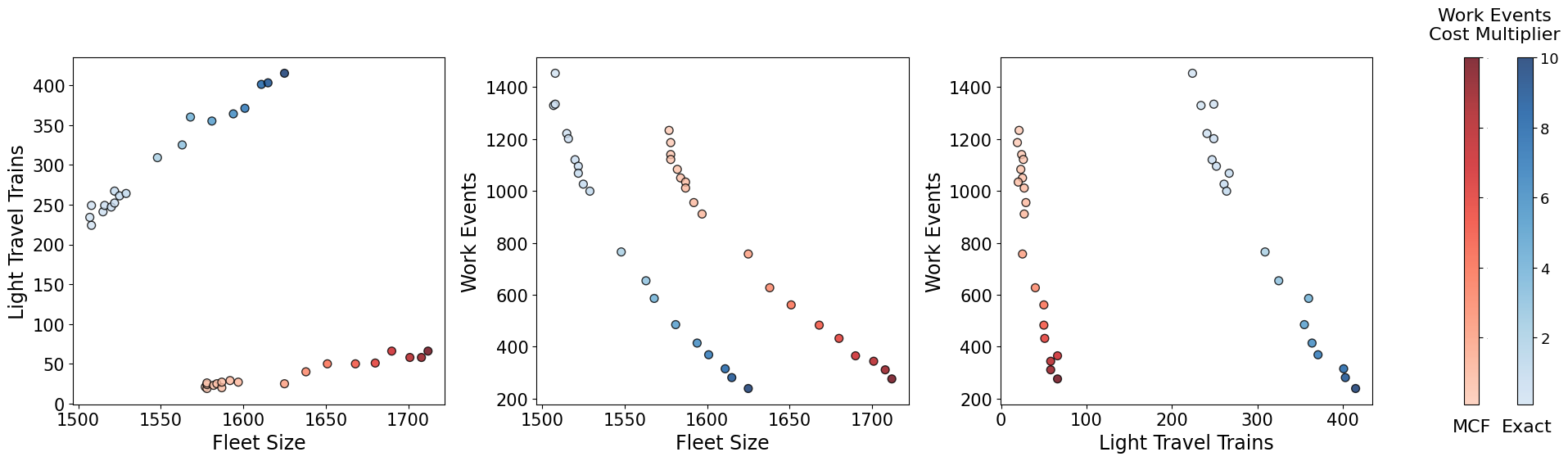}} \\
    \subfigure[Varying Relocation Cost ($g$) \label{fig:sens_reloc}]{\includegraphics[width=\textwidth]{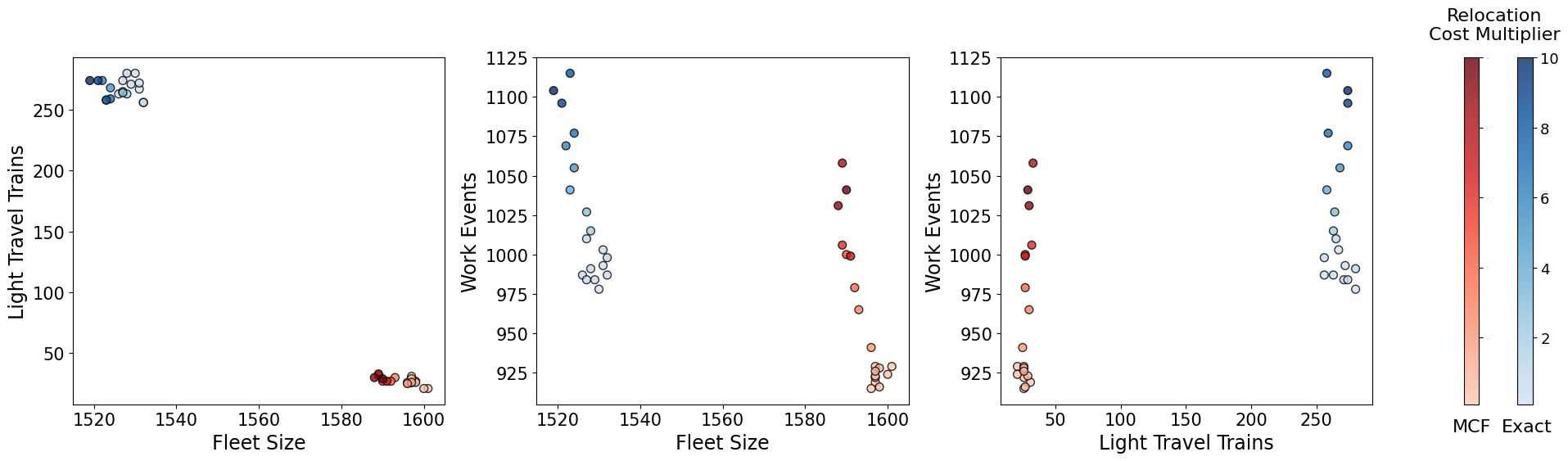}}
    \caption{Impact of Varying Cost Coefficients on Light Travel, Fleet, and Work Events}
    \label{fig:senstanalys}
\end{figure}

The four key cost components examined are: fleet size, work events, light traveling, and repositioning of inactive units. For each component, its associated cost coefficient is independently scaled from 0.1 to 1.0 (step size 0.1) and from 1.0 to 10.0 (step size 1.0). Figure \ref{fig:senstanalys} highlights the resulting trade-offs across three key performance indicators: fleet size, number of work events, and number of light travel arcs. 

\vspace{-0.5cm}
\subsubsection{Fleet Sizing vs. Locomotive Utilization}
Fleet cost ($q$) dominates the default cost configuration and strongly shapes the solution structure. While minimizing capital investment favors a leaner fleet, this must be weighed against the operational burden of over-utilization. This reflects the fundamental trade-off between \textit{capacity}--owning more locomotives--and \textit{efficiency}--maximizing each unit's use. 

As shown in Figure \ref{fig:sens_fleet}, when $q$ is low, the model prefers to expand the fleet and reduces reliance on repositioning, whether via work events or light travel. Similar behavior emerges when work event costs ($c$) or crew cost ($e$) is high as shown in Figures~\ref{fig:sens_lt} and~\ref{fig:sens_wc}: relocation becomes expensive, making it advantageous to stage locomotives closer to demand. On the other hand, as $q$ increases, the fleet contracts and repositioning becomes essential to meet service requirements. The same pattern arises when $e$ or $c$ are low, incentivizing aggressive repositioning in lieu of capital expansion.

\begin{figure}[H]
    \centering    \includegraphics[width=0.5\linewidth]{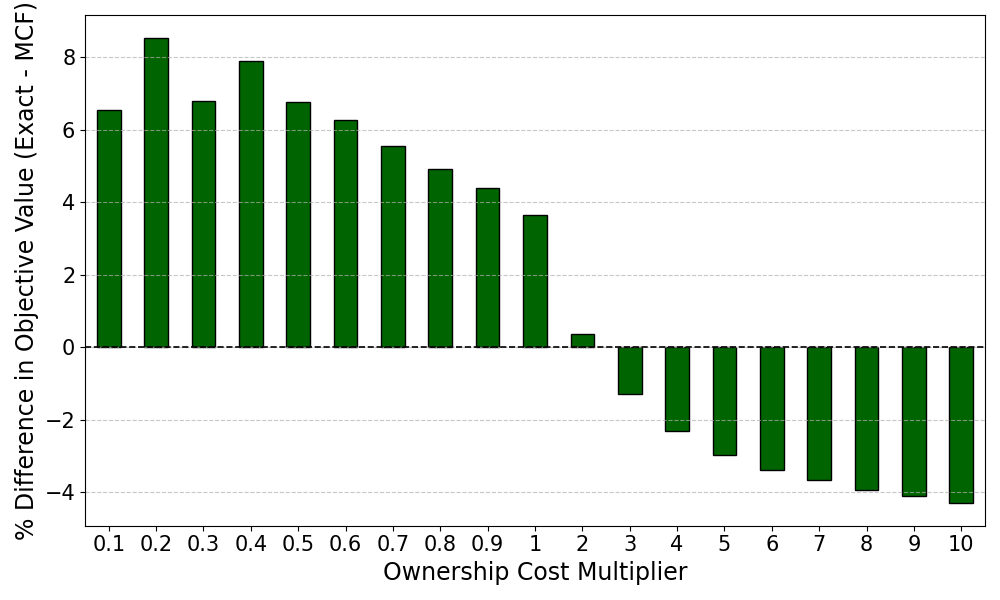}
    \caption{Objective Value Difference of Exact Method from MCF (\%) under Varying Ownership Cost ($q$)}
    \label{fig:FC_proportional_objs}
    \vspace{-0.3cm}
\end{figure}

Between the two models, the exact optimization demonstrates greater sensitivity to cost parameters and consistently achieves smaller fleet sizes by up to 90 fewer units across scenarios. 
This capability is especially valuable in capacity-constrained settings, such as during seasonal demand peaks with limited inventory. 
Under high-$q$ scenarios, the exact model begins to outperform MCF in total objective value starting at cost factor 3, with up to 4\% savings at cost factor 10 (Figure \ref{fig:FC_proportional_objs}). This advantage stems from its ability to dynamically leverage a richer set of light travel arcs, whereas the MCF model plateaus in performance due to its rigid $A_L$ structure.

\vspace{-0.5cm}
\subsubsection{Light Travel Optionality vs. Network Saturation}
Light travel arcs offer critical flexibility in locomotive repositioning, especially under uncertainty. While DH is more cost-effective, it relies heavily on the timely execution of scheduled train moves, thus assuming a reliable sequence for repositioning. In reality, disruptions, such as infrastructure failures or weather, can sever these chains. A single missed arrival may cascade into downstream locomotive shortages and trigger widespread delays. Light travel mitigates this fragility by allowing direct repositioning independent of train schedules. This optionality enhances system resilience and decentralization--particularly valuable in volatile environments. 
However, light travel is costly. In excess, it contributes to congestion by adding to traffic density and scheduling conflicts. The central trade-off, then, lies between \textit{repositioning flexibility} and \textit{network saturation}.

Crew cost, the dominant component of light travel cost ($e$), directly governs arc usage. As shown in Figure \ref{fig:sens_lt}, light travel declines sharply as $e$ increases. At high $e$, the exact model converges to a minimal arc set, comparable to that of the MCF model. At low $e$, however, the exact method activates over 600 arcs, while MCF maintains its limited, predictable arc usage. Although MCF benefits from faster solve times and quality solutions thanks to its compact, narrowly filtered arc set,  its inflexibility prevents it from capitalizing on low $e$.
By contrast, the exact model dynamically adjusts arc usage in response to shifting conditions. 
At high $e$, it emulates MCF behavior, demonstrating that flexibility is selectively invoked when cost-justified.  

\begin{figure}[!]
  \centering
  \subfigure[$A_L$ Generated with MCF method]{%
    \includegraphics[width=0.48\linewidth]{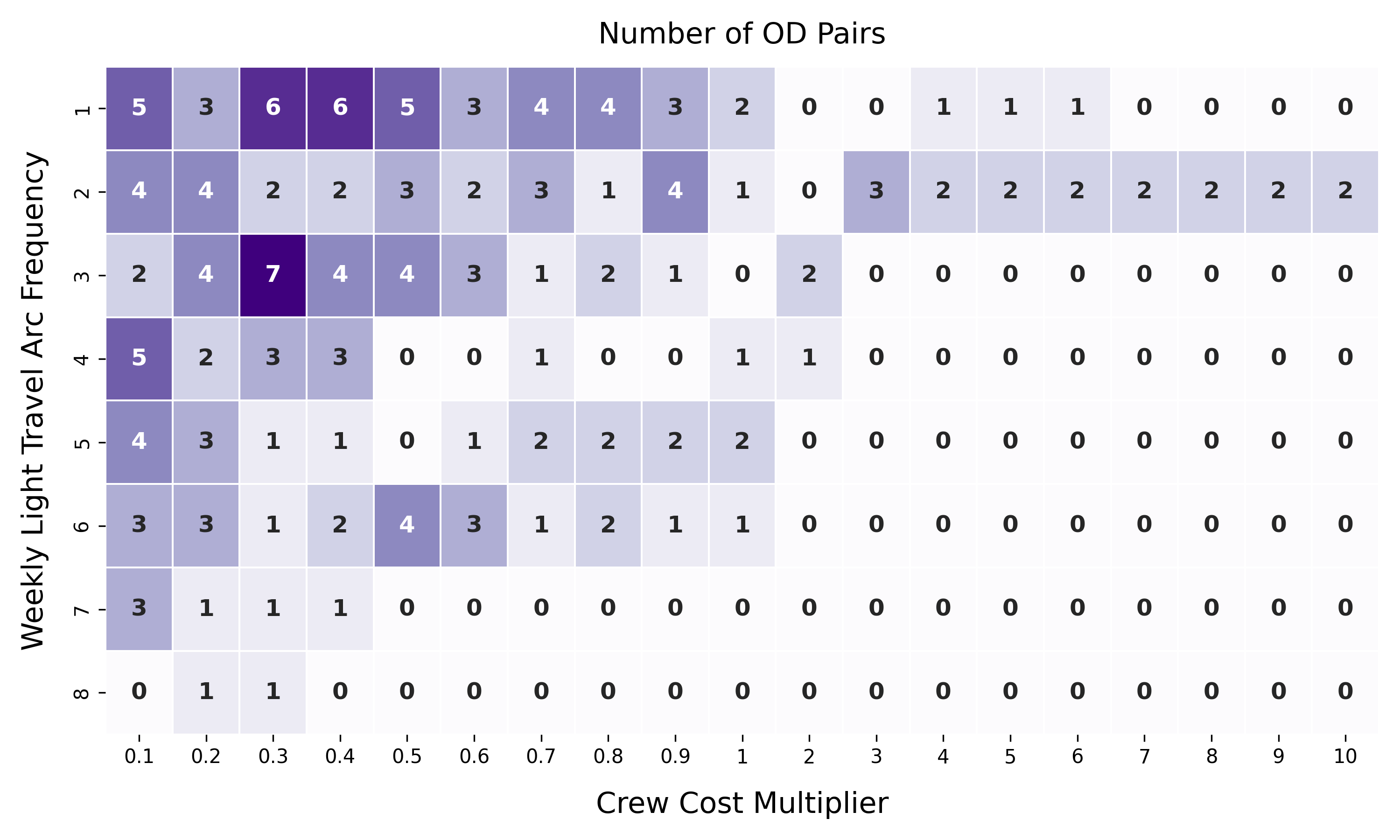}
    \label{fig:cc_exact_od}
  }
  \hfill
  \subfigure[$A_L$ Generated with Exact method]{%
    \includegraphics[width=0.48\linewidth]{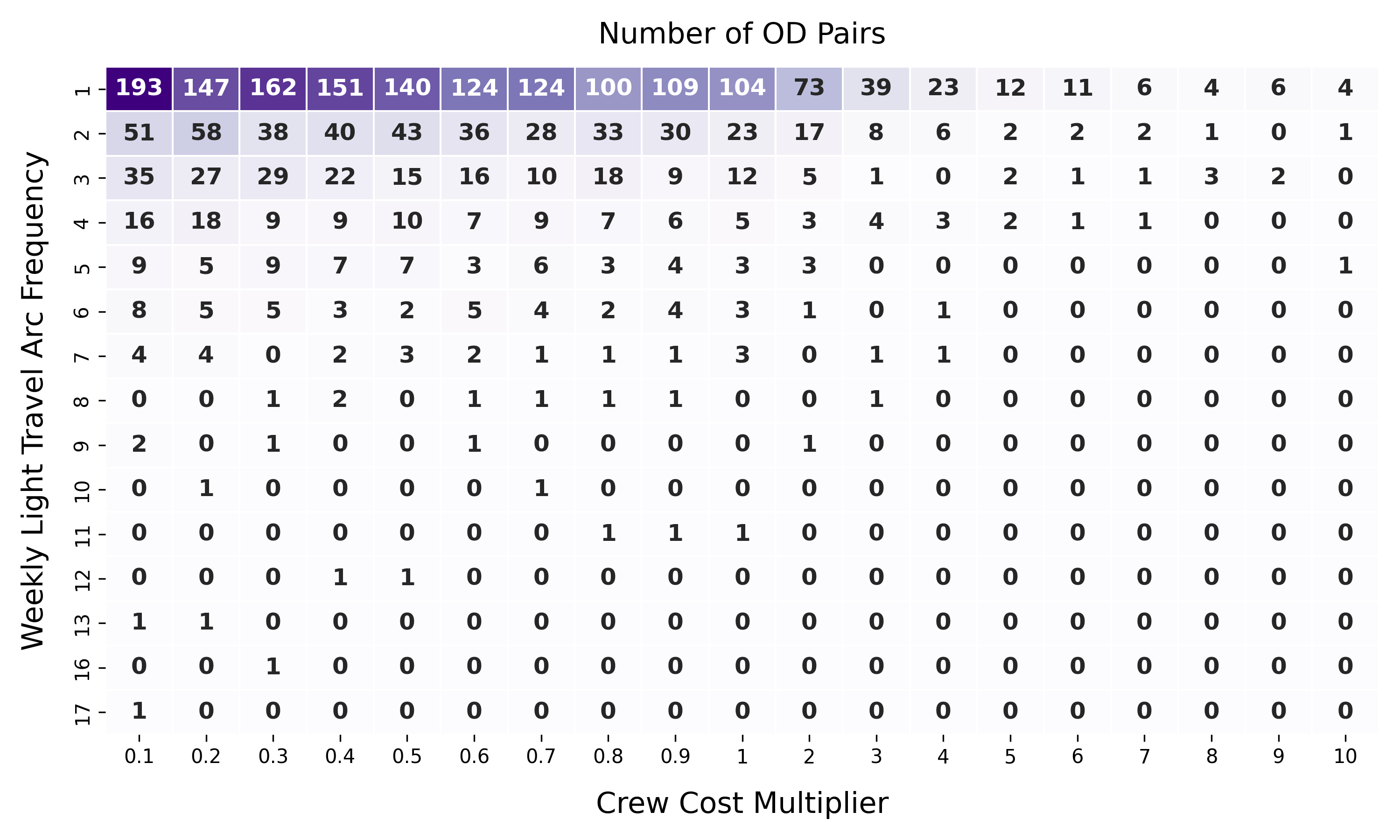}
    \label{fig:cc_MCF_od}
  }
  \caption{Number of Unique Origin-Destinations Traveled by Light Travel Arcs Under Varying Crew Cost ($e$)}
  \label{fig:cc_od_frequency}
\end{figure}

\begin{figure}[
!]
  \centering
  \includegraphics[width=0.5\linewidth]{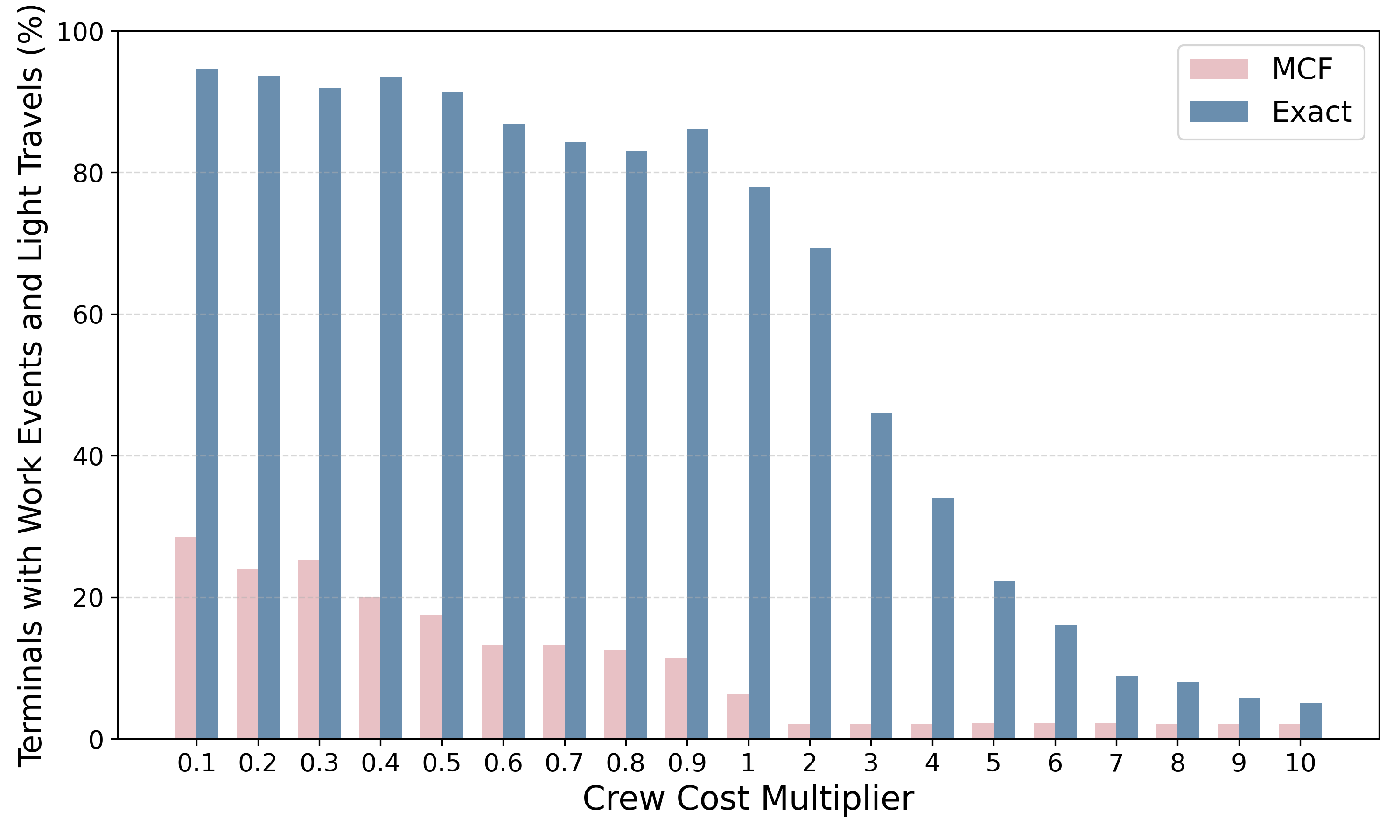}
    \vspace{-0.3cm}
    \caption{Percentage of Terminals with Concurrent Work Events and Light Travels}
    \label{fig:cc_coverage}
    \vspace{-0.3cm}
\end{figure}

Figure \ref{fig:cc_od_frequency} displays how light travel frequency between unique OD terminal pairs varies with $e$. At low $e$, the exact model activates a broad set of OD pairs at low frequency; as $e$ increases, arc usage becomes more concentrated and MCF-like. OD pairs that persist at high $e$ mark critical relocation corridors. Flexible light travel also enables broader coordination with work events. As Figure~\ref{fig:cc_coverage} shows, low $e$ increases the overlap between terminals that have work events and also serve as light travel endpoints; this is an indication of decentralized repositioning, where a wide range of terminals participate in repositioning. The exact model achieves far greater overlap, unlocking terminals MCF never uses. MCF's light travels are confined to a select surplus-deficit pairs, limiting adaptability. 
As $e$ rises, this overlap declines, indicating a shift toward a more concentrated repositioning pattern that relies on a small subset of key terminals. The system becomes less willing to execute multi-hop repositioning that couples work events with light travel and shifts burden to work events alone.

\subsubsection{Operational Overhead vs. Frequent Repositioning} 
A core motivation of this study is to minimize unproductive locomotive movement. Keeping all units attached for an entire train journey often leads to over-provisioning and idle power (under-utilization). Decentralized repositioning--dropping off or picking up at intermediate stops--enables timely reuse and reduces DH. However, these shorter repositioning activities come at the expense of increased operational overhead: more terminal activity, tighter scheduling, and crew coordination. The resulting tradeoff is \textit{repositioning frequency} and \textit{operational complexity}. 

\begin{figure}[!]
    \centering
    \subfigure[Active Terminals vs. Coverage ($c$)]{
        \includegraphics[width=0.29\textwidth]{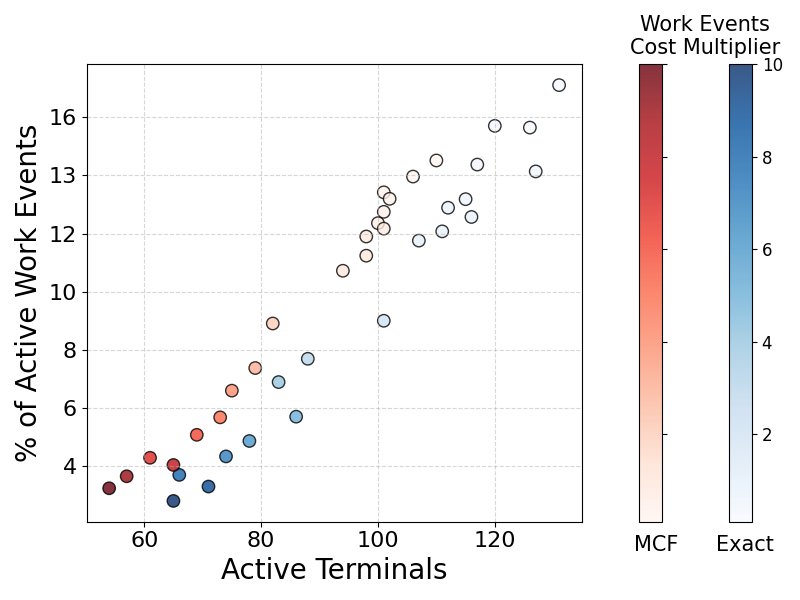}
        \label{fig:WC_coverageplot}
    }
    \hfill
    \subfigure[Number of Work Events and Locomotives Involved ($c$)]{
        \includegraphics[width=0.65\textwidth]{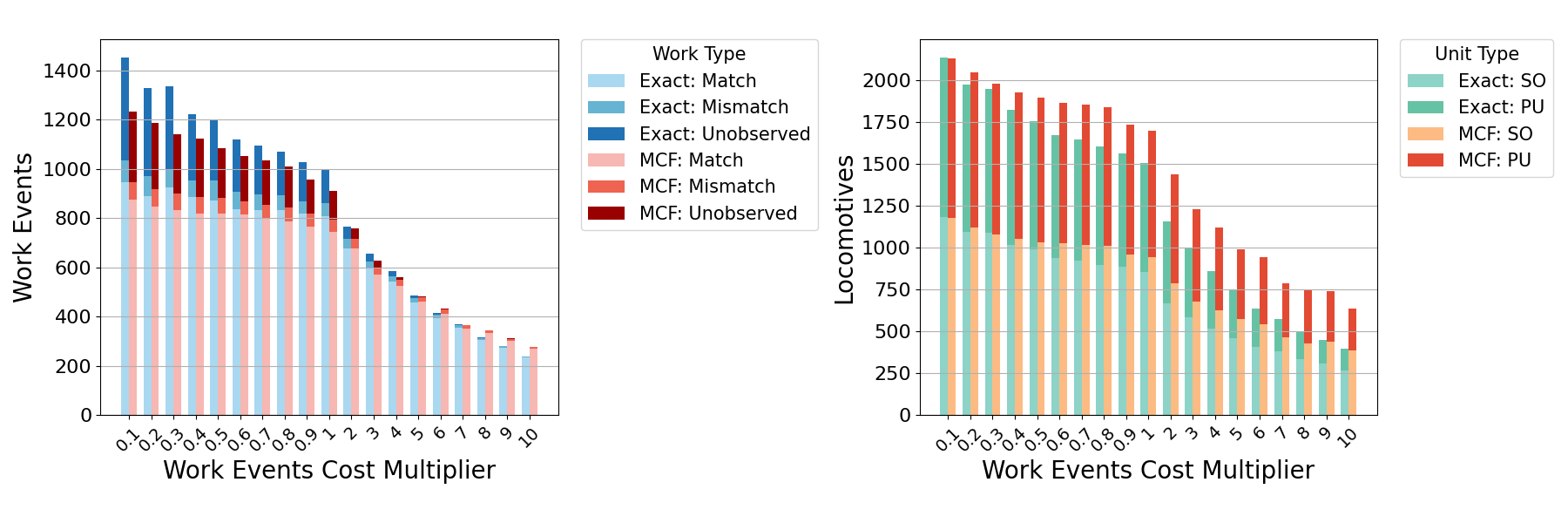}
        \label{fig:wc_events_units}
    }
    
    \vspace{0.1cm}
    
    \subfigure[Active Terminals vs. Coverage ($g$)]{
        \includegraphics[width=0.29\textwidth]{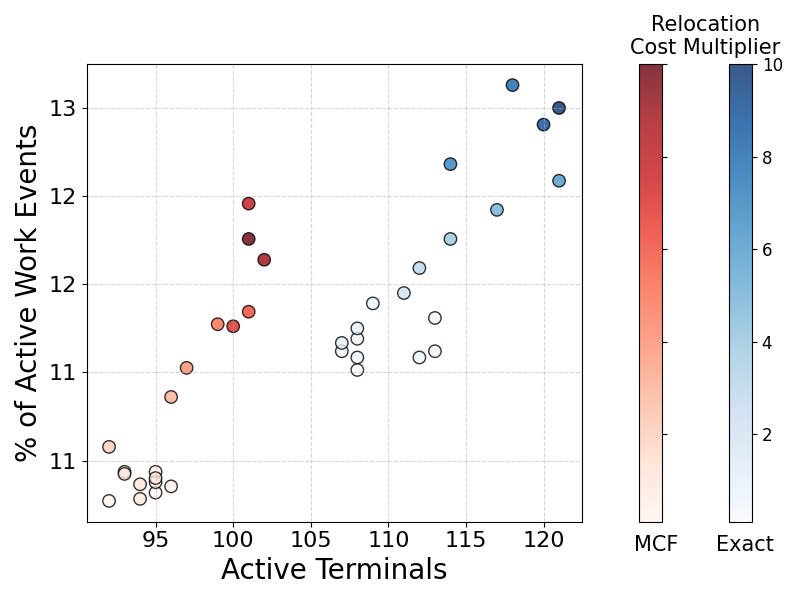}
        \label{fig:RC_coverageplot}
    }
    \hfill
    \subfigure[Number of Work-Events and Locomotives Involved ($g$)]{
        \includegraphics[width=0.65\textwidth]{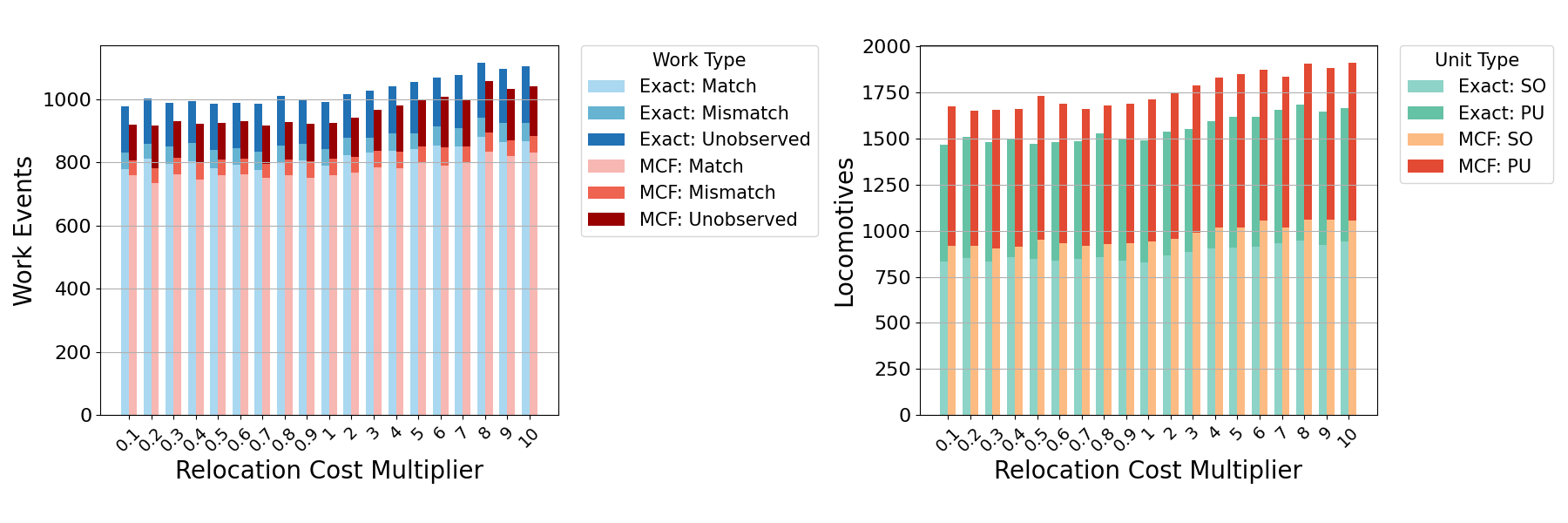}
        \label{fig:rc_events_units}
    }

    \caption{Number of Terminals and Locomotives undergoing Work Events along with Coverage Ratio}
    \label{fig:WC_RC_figures}
\end{figure}

Figure \ref{fig:WC_coverageplot} illustrates this clearly as work event costs $c$ are varied. When $c$ is low, the system favors decentralized repositioning over full-route commitments. This results in higher number of terminals active with work events and work event coverage ratio, which measures the portion of utilized work events out of total possible event opportunities. As $c$ increases,
total work events and participating units decline (Figure~\ref{fig:wc_events_units}). Interestingly, MCF consistently involves more locomotives across fewer or slightly more events, whereas the exact model has smaller number of units across a broader spread of repositioning points.
A similar pattern emerges when scaling relocation cost $g$ (\ref{fig:RC_coverageplot},  \ref{fig:rc_events_units}): the exact model consistently executes more events with fewer units, reinforcing its decentralized repositioning strategy. Eventually, the number of active terminals plateaus while coverage continues rising, suggesting denser work event frequency at already-active terminals rather than activating new ones. 

\begin{figure}[h]
  \centering
  \subfigure[Varying Work Events Costs ($c$)]{%
    \includegraphics[trim=0 0 220 0, clip, width=0.415\linewidth]{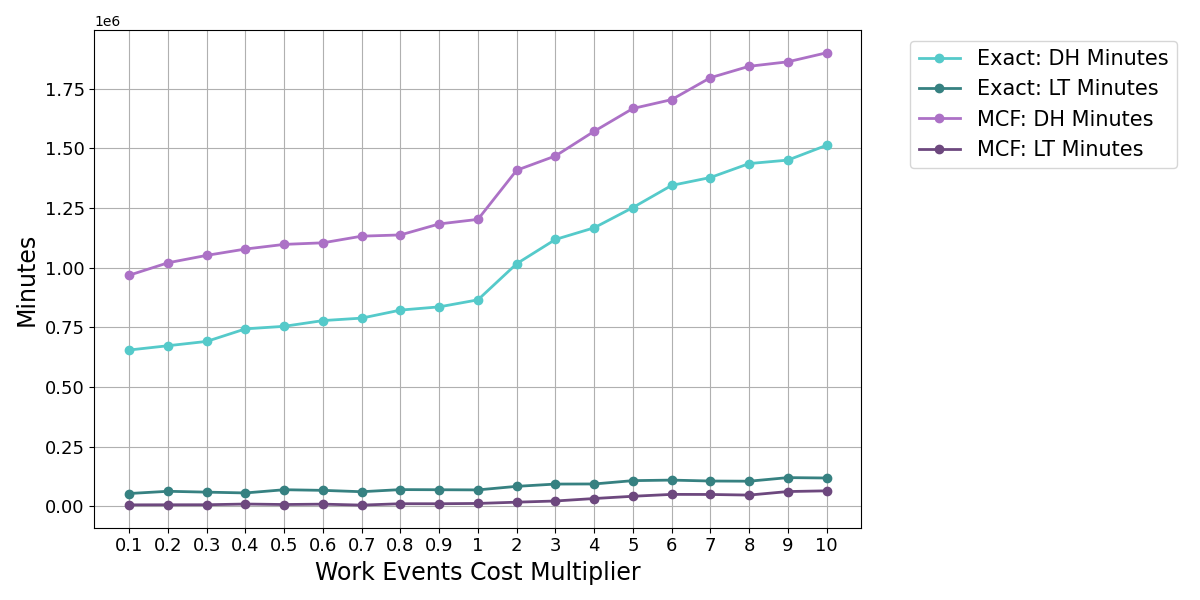}
    \label{fig:wc_dg_mins}
  }
  \hfill
  \subfigure[Varying Relocation Costs ($g$)]{%
    \includegraphics[width=0.555\linewidth]{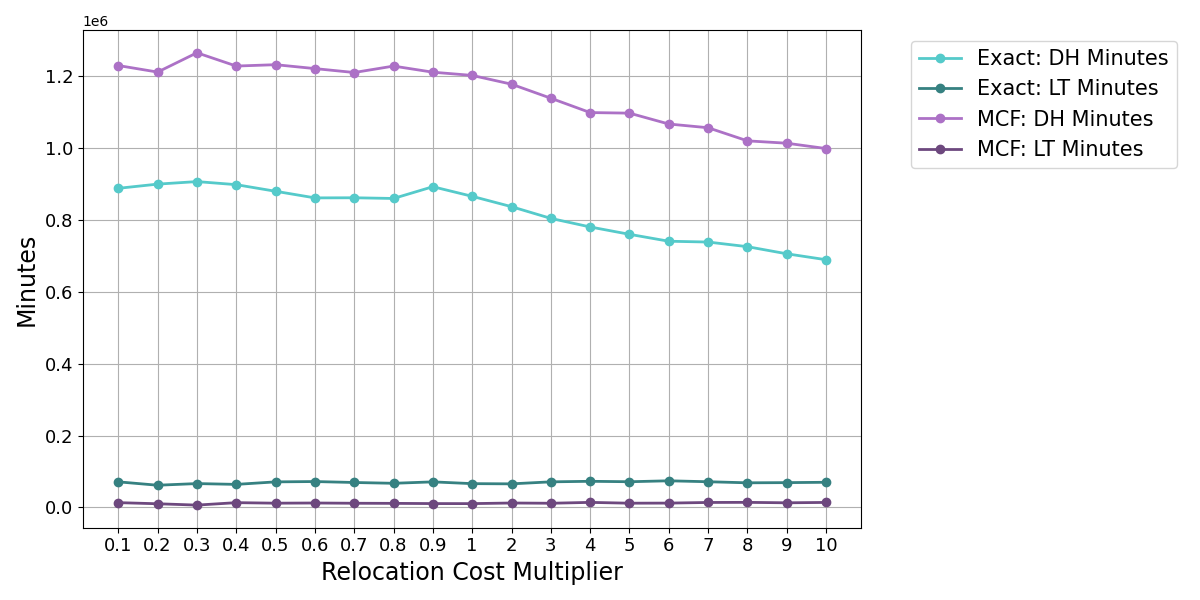}
    \label{fig:rc_dg_mins}
  }
  \caption{Repositioning Minutes by Mode}
  \label{fig:DG_mins}
  \vspace{-0.2cm}
\end{figure}

Figure \ref{fig:DG_mins} further complements this trend: as $c$ increases, DH minutes rise sharply from short-range repositioning, while light travel minutes grow only modestly. When $g$ increases, DH declines but light travel remains flat, hinting diminishing returns from controlling light travel. Thus, decentralized repositioning reduces DH by leveraging network flexibility but introduces a substantial operational strain. 

\subsection{Impacts of Practice-Based Constraints on System Performance}

All previous results assume full flexibility in scheduling work events (V0), representing a theoretical upper bound. In practice, railroad industries face structural limits such as fixed infrastructure and policy constraints. To reflect these realities, this section evaluates system performance under a series of operational constraints introduced in Section \ref{sec:3.3}. All experiments limit terminals to at most 6 work events per day ($\theta$).

\subsubsection{Capacity Extension in Baseline Plan}
This setting adopts the most conservative approach: optimization is restricted to only those terminal-days that had work events in the baseline plan. The model can increase the number of events at these terminal-days, up to the minimum of $h_{k,d} + \lambda$ (original allocation plus extra capacity) or $\theta$.

\begin{figure}[htbp!]
    \vspace{-0.2cm}
    \centering

    \subfigure[$\text{V1}_1 (\lambda=1)$]{
        \includegraphics[width=0.31\linewidth]{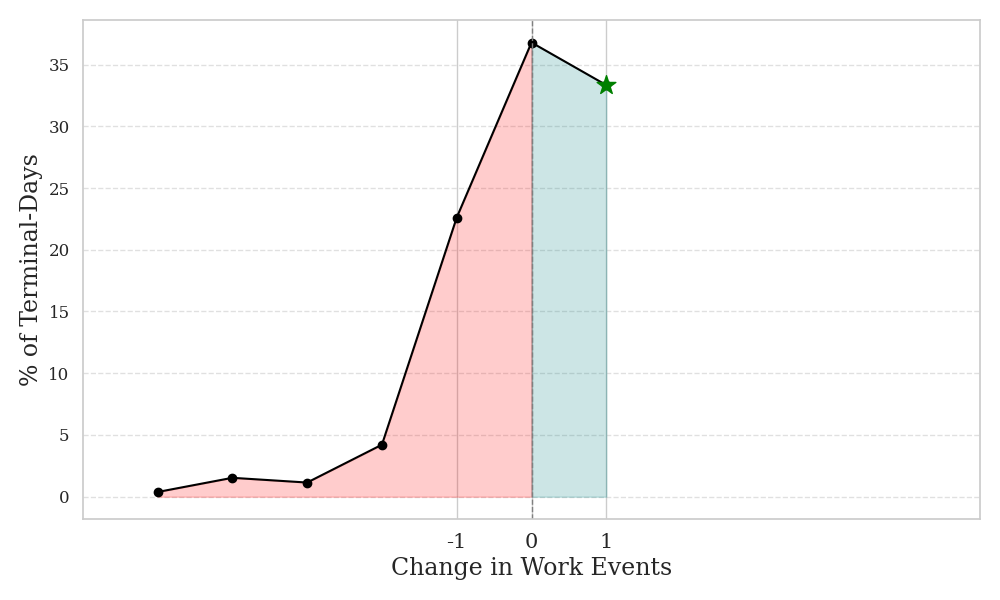}
    }
    \hfill
    \subfigure[$\text{V1}_2 (\lambda = 2)$]{
        \includegraphics[width=0.31\linewidth]{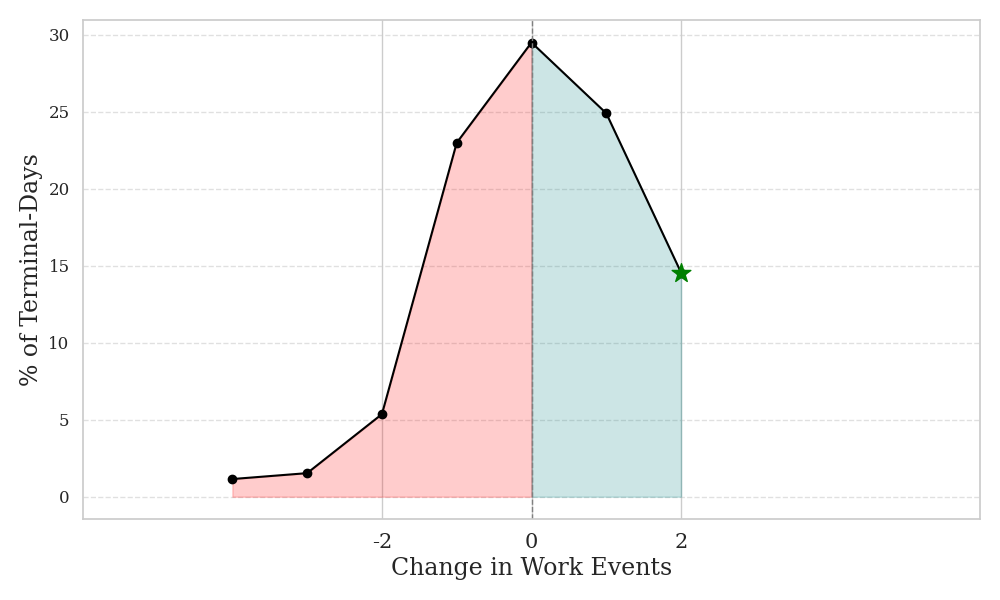}
    }
    \hfill
    \subfigure[$\text{V1}_3 (\lambda = 3)$]{
        \includegraphics[width=0.31\linewidth]{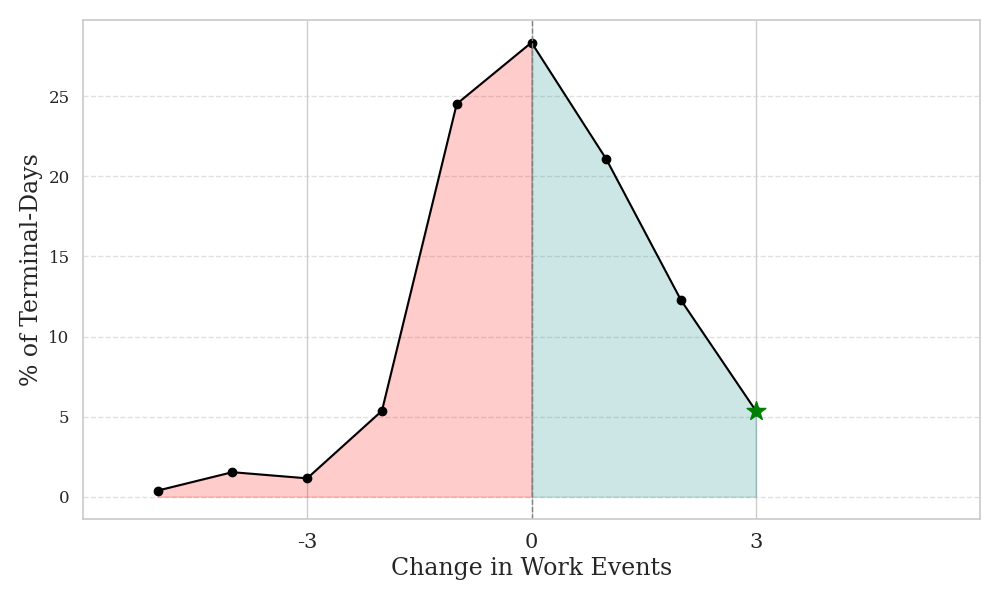}
    }

    \makebox[\textwidth][c]{%
    \subfigure[$\text{V1}_4 (\lambda = 4)$]{%
        \includegraphics[width=0.31\linewidth]{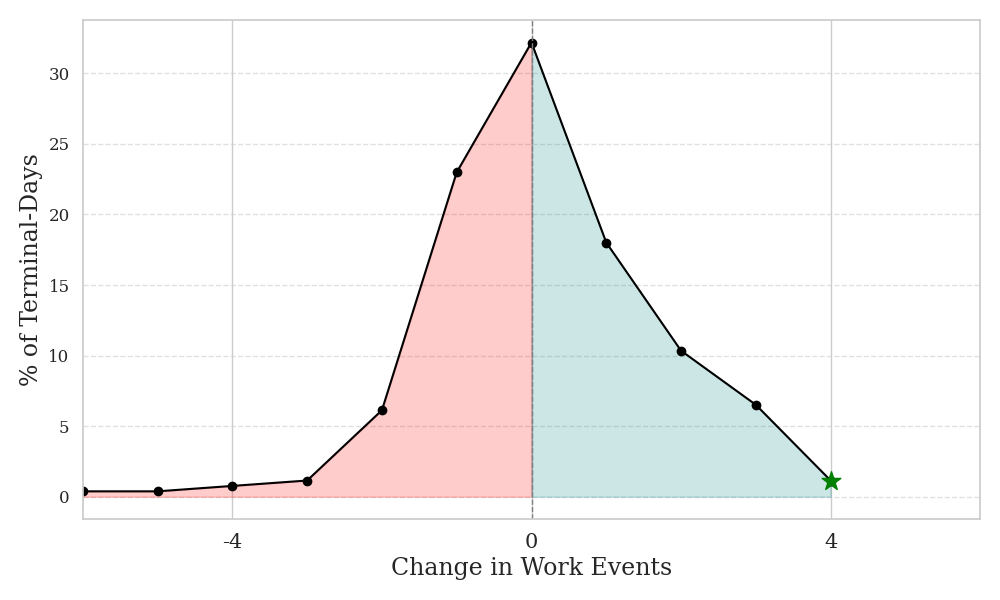}
    }
    \hspace{0.06\linewidth}
    \subfigure[$\text{V1}_5 (\lambda = 5)$]{%
        \includegraphics[width=0.31\linewidth]{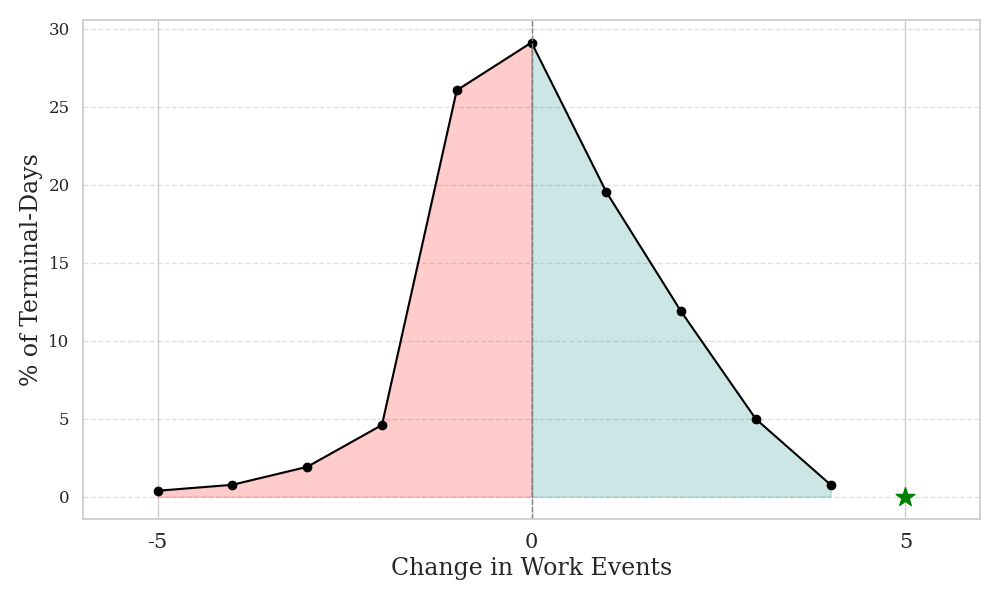}
    }
}
    \caption{Distribution of Work Event Changes (LAP – Baseline) across Baseline-Active Terminal-Days}
    \label{fig:V1_util_improve}
    \vspace{-0.5cm}
\end{figure}

Figure \ref{fig:V1_util_improve} shows how terminal-day's capacity utilization shifts with increasing $\lambda$. Red shading indicates terminal-days where the optimized number of work events falls below or matches the original level ($h_{k,d}$); they underutilize even the original assignment. Blue shading reflects added activity, and star markers show terminal-days that reach their full capacity limit. 
As $\lambda$ grows, fewer terminal-days hit that cap, and by $\lambda = 5$, none do. This indicates expanding capacity at already-active sites yield only diminishing returns.

When examining objective improvements, marginal returns diminish after $\lambda=2$, implying most optimization potential is captured with modest extension (more details in Appendix). This motivates exploring strategies that activate new terminals or terminal-days.
For the remaining experiments, $\text{V1}'$ is adopted as the baseline plan. It doubles the allocation $h_{k,d}$ for active terminal-days, capped at $\theta$. $\text{V1}'$ represents a pragmatic operational stretch and serves as a warm-start point for the activation-based variants (V2–V5). This level was specified by industry partners as the most feasible near-term expansion.

\subsubsection{Incremental Augmentation vs. Network-wide Redesign}
To assess how flexibility in work event scheduling affects system efficiency, two strategies are evaluated. 
The \emph{incremental} approach (V2, V3) adds a limited number of terminals ($\alpha^C$) or terminal-days ($\alpha^D$) on top of the baseline ($\text{V1}'$). The \emph{network-wide redesign} approach (V4, V5) selects terminals ($\alpha^E$) and terminal-days ($\alpha^F$) from scratch, unconstrained by the baseline layout. These $\alpha$-thresholds across the two strategies match the total number of active elements (e.g, $\alpha^E$ = baseline-active terminals + $\alpha^C$). All models are warm-started from $\text{V1}'$ solution for consistency and computational efficiency, though this may introduce bias toward baseline structures. Additionally, higher-$\alpha$ models are initialized from the preceding lower-$\alpha$ solutions to facilitate marginal comparisons.

\begin{figure}[H]
    \centering   
    \subfigure[Objective Improvements under Limiting Number of Active Terminals]{
        \includegraphics[width=0.8\linewidth]{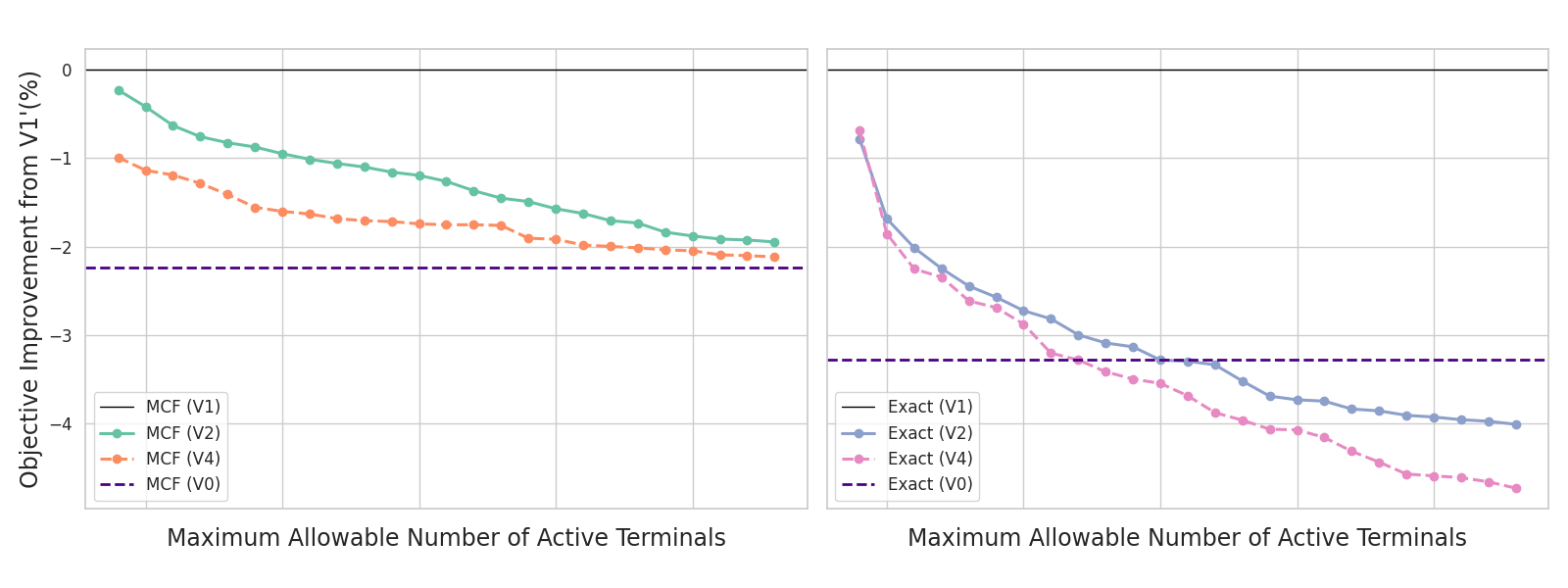}
        \label{fig:CE_obj}
        }
    
    \subfigure[Objective Improvements under Limiting Number of Active Terminal-Days]{
    \includegraphics[width=0.8\linewidth]{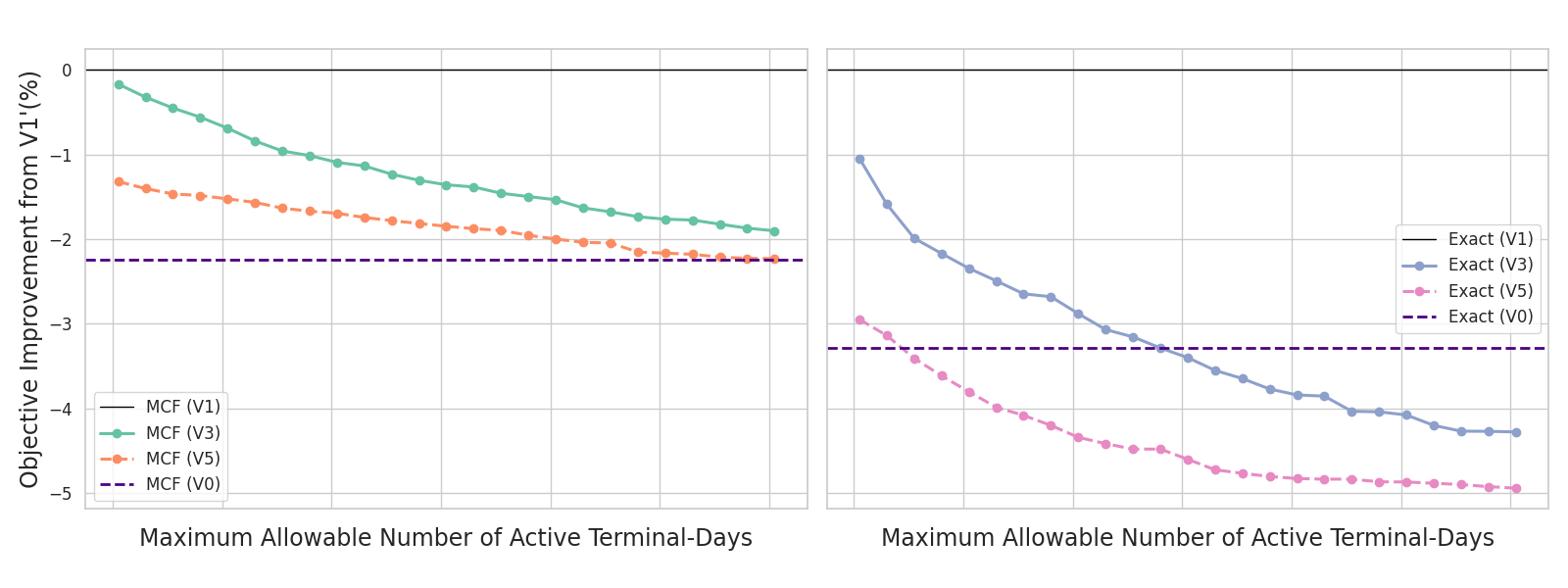}
    \label{fig:DF_obj}
    }
    \caption{Objective Improvements of Model Variants (V2, V3, V4, V5) relative to $\text{V1}'$}
    \label{fig:cd_ef_obj}
    \vspace{-0.2cm}
\end{figure}

As illustrated in Figures \ref{fig:CE_obj} and \ref{fig:DF_obj}, the redesign models (V4, V5) consistently outperform their incremental counterparts (V2, V3) in their objective value across all tested activation thresholds ($\alpha$). Improvements taper off as $\alpha$ increases, indicating diminishing returns. 
Redesigns yield stronger solutions earlier by constructing globally coherent configurations from the outset. On the other hand, incremental models must work around the baseline's inefficiencies, and continually layer improvements. 
All variants eventually outperform the unconstrained V0 solution (solved for 24-hours without warm-starting or work event constraints), due to cumulative refinement and guided intializations.

The advantage of redesign approach does not come merely from replacing underperforming components in isolation, but from coordinated substitutions across the network and their interaction.
Incremental strategies are locked into spatial and temporal structure of baseline work events, which limit their ability to integrate high-value additions.
Redesigns can include individually suboptimal components if they contribute to a more effective configuration when optimized jointly.

These findings underscore a strategic tradeoff. Incremental models offer continuity from baseline with lower disruption but their benefits plateau. Redesigns require more computational effort and structural change, but deliver better outcomes. 
A hybrid policy may balance practicality and performance: begin with incremental steps, then selectively reallocate baseline components.

\subsubsection{Spatial vs. Spatio-Temporal Flexibility}
Another key distinction across variants lies in the level of control over when and where work event can happen. \emph{Spatial} strategies (V2, V4) activate terminals and it is active for all days if used at all, whereas \emph{spatio-temporal} models (V3, V5) activates terminals for specific days.
This finer control allows to target terminals that are high-impact on specific days, capturing intermittent opportunities that spatial models overlook as they prioritize sustained, week-long usage.

\begin{figure}[htbp]
    \vspace{-0.5cm}
    \centering
    \subfigure{%
        \includegraphics[width=0.47\textwidth]{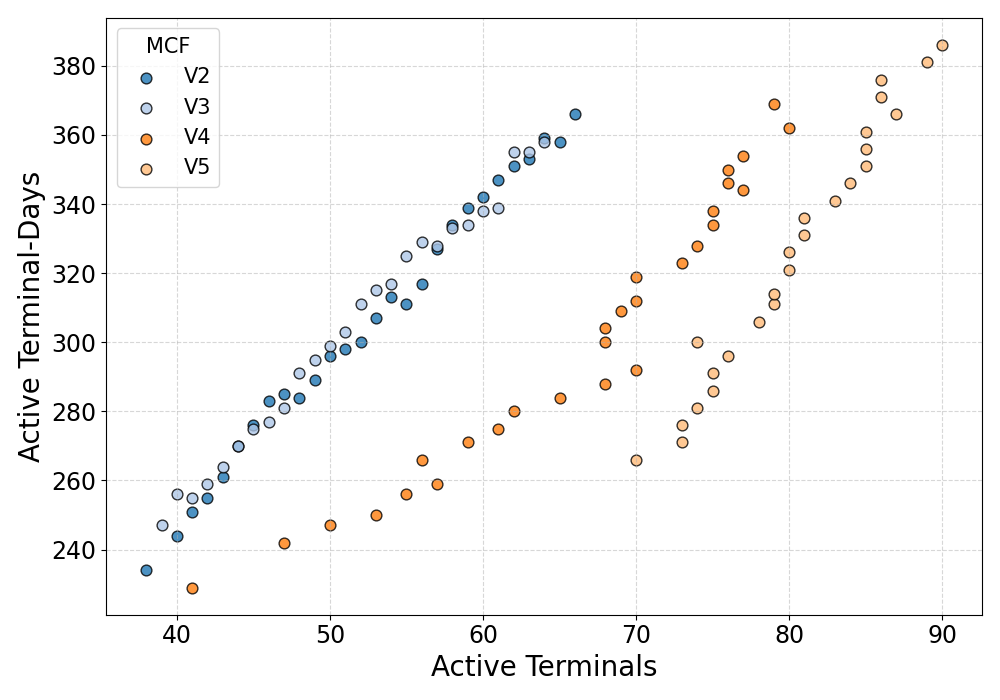}
        \label{fig:2Dscatter_1}
    }
    \hfill
    \subfigure{%
        \includegraphics[width=0.47\textwidth]{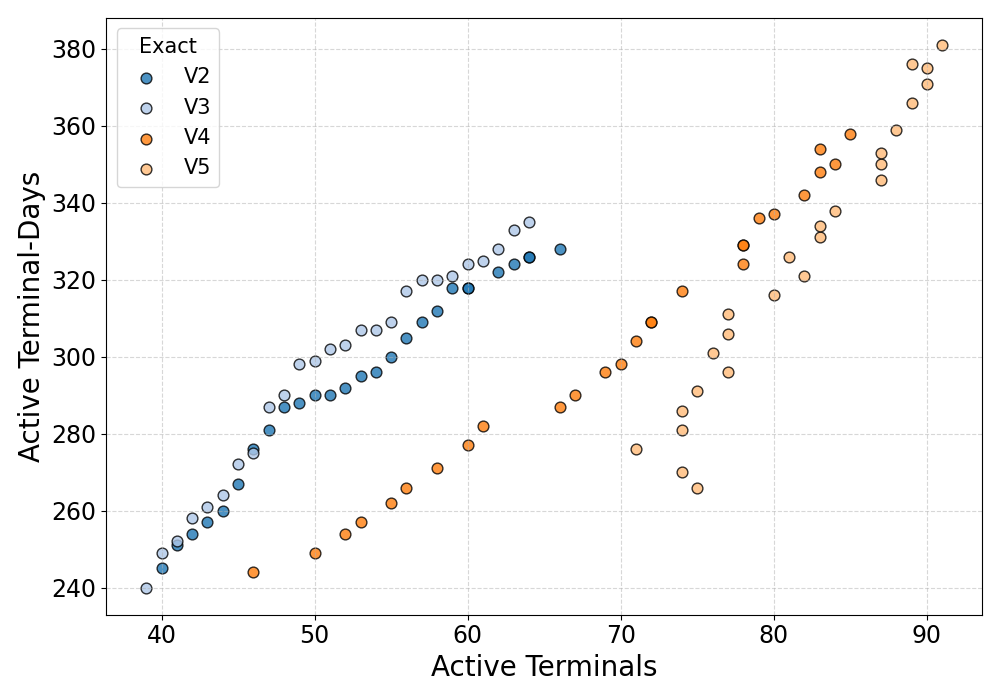}
        \label{fig:2Dscatter_2}
    }
    \caption{Number of Active Terminals and Terminal-Days Across Model Variants}
    \label{fig:V2-5:terms_tds}
    \vspace{-0.2cm}
\end{figure}

Spatio-temporal models achieve broader spatial reach and better temporal coverage. 
As shown in Figure \ref{fig:V2-5:terms_tds}, they tend to engage more unique terminals than their spatial counterparts, even when total active terminal-day counts are similar. 
They also make fuller use of their allowed terminal-day budgets by activating more terminal-days when the same amount is allowed, suggesting more efficient targeting of high-value time windows.

\begin{figure}[H]
    \centering
    \includegraphics[width=0.7\linewidth]{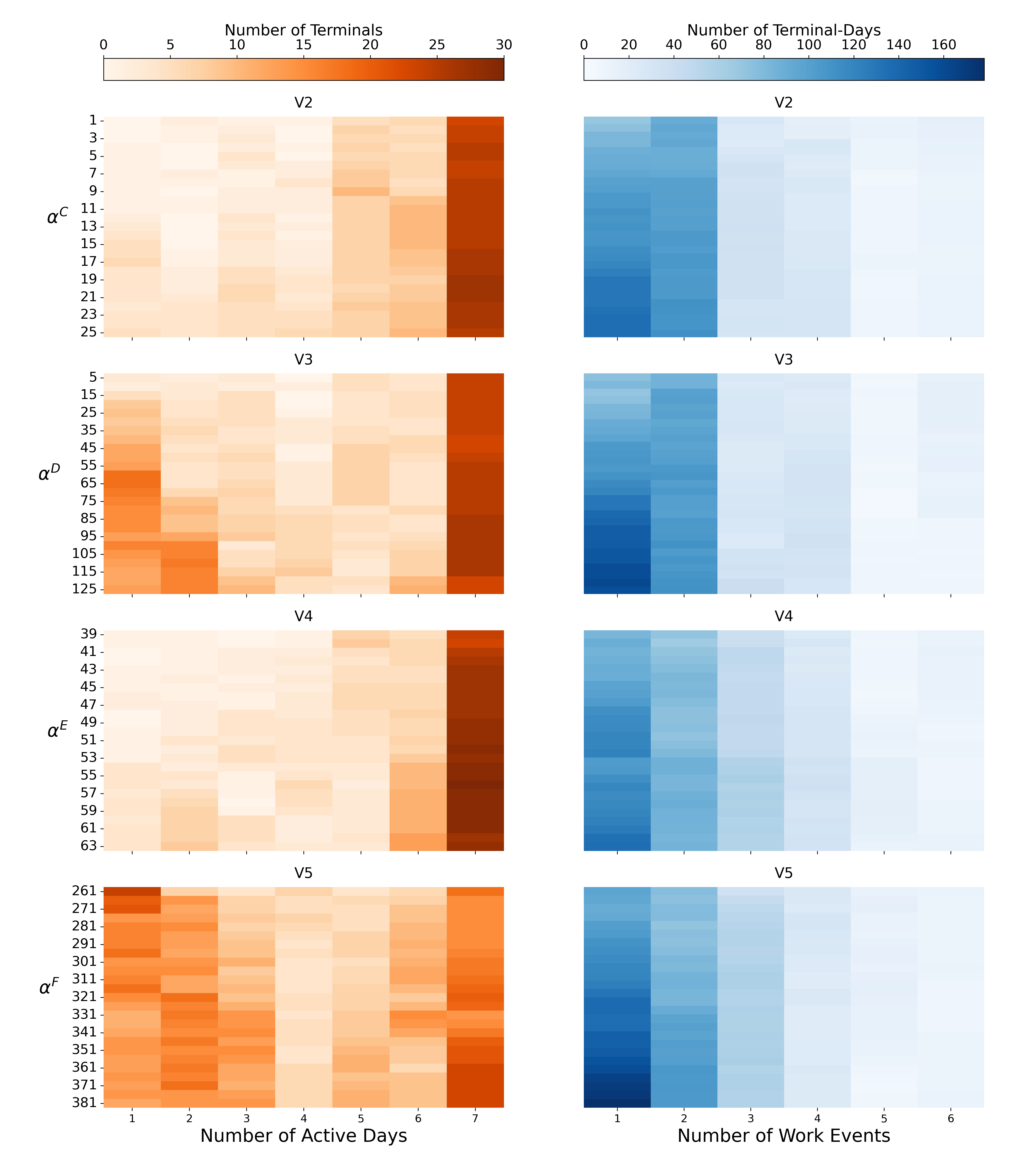}
    \caption{Spatio-Temporal Activation Patterns Across Model Variants Under Exact Model}
    \label{fig:exact_heatmaps}
    \vspace{-0.4cm}
\end{figure}

Load distribution patterns differ as well. Spatial models concentrate work events at a group of high-performing terminals, as visible in the left-hand panels of heatmaps of Figure \ref{fig:exact_heatmaps}. V2 and V4 display vertical bands for 7-days, a sign of full-week usage. Spatio-temporal models (especially V5) exhibit more diffuse activation, spreading work events across both terminals and days. This diversification helps reduce overload at any single terminal and contributes to greater cost savings. The corresponding heatmap for the MCF method has a similar pattern and is provided in the Appendix for reference.

At higher $\alpha$-thresholds, spatio-temporal models often underutilize their terminal-day budgets ($\alpha^D$ or $\alpha^F$), which indicates that most gains are realized from a small set of high-value space-time combinations.
In contrast, spatial models tend to use their full allocation ($\alpha^C$ or $\alpha^E$), extending to those with marginal benefit once broader coverage is exhausted.

\section{Conclusion} \label{sec:5}
This study presents a strategic-level formulation of the Locomotive Assignment Problem that directly integrates the selection of work event locations and incorporates a broad set of operational constraints. In doing so, it takes a meaningful step toward bridging the gap between optimization modeling and real-world application. 
The study proposes tractable reduction methods that allow exact optimization over an extensive set of light travel options. Although the exact method is still computationally intensive, it establishes a ground truth for evaluating solution quality. It provides valuable structural insights and serves as a true benchmark for evaluating the heuristic’s effectiveness. The MCF heuristic is validated to perform well, offering high-quality solutions with tractable solve times. Together, these results offer both a rigorous analytical foundation and practical value to decision-makers, to navigate the tension between locomotive capacity, repositioning flexibility, and operational overhead.
Furthermore, it is demonstrated that the exact method can outperform MCF when locomotive ownership cost is high, e.g., for highly-constrained systems.

Future work can build on this strategic framework by integrating its optimized outputs into the subsequent tactical and operational planning phases as input parameters. There is also considerable opportunities to refine the light travel arc selections through decomposition techniques. In parallel, modeling the routing and chaining of light travel arcs into feasible locomotive journeys can better capture the temporal dynamics of repositioning. Finally, introducing operational uncertainty would allow this framework to evolve into a resilient decision-support system under real-world volatility.

\if0\blind{
\section*{Acknowledgements}
This research was partly funded by NSF AI4OPT under Award 2112533. } \fi


\bibliographystyle{apalike}
\spacingset{1}
\bibliography{citation.bib}

\spacingset{1.5}
\appendix

\newpage
\section*{Appendix} \label{Appendix}
\section{Generating Exact Set of Light Travel Arcs} \label{Appendix:6.1}

To reduce the full set of candidate light travel arcs, the following reduction steps may be applied in order:
\begin{enumerate}
    \item Earliest Reachability: Each arrival-ground node is connected to the first available ground-departure node at each destination terminal that it can reach.
    \item Latest Origin Filtering: If multiple arrival-ground nodes from the same origin terminal connect to the same ground-departure node, only the arc from the latest such origin node is retained.
\end{enumerate}

\begin{figure}[!]
    \centering
    \includegraphics[width=0.8\linewidth]{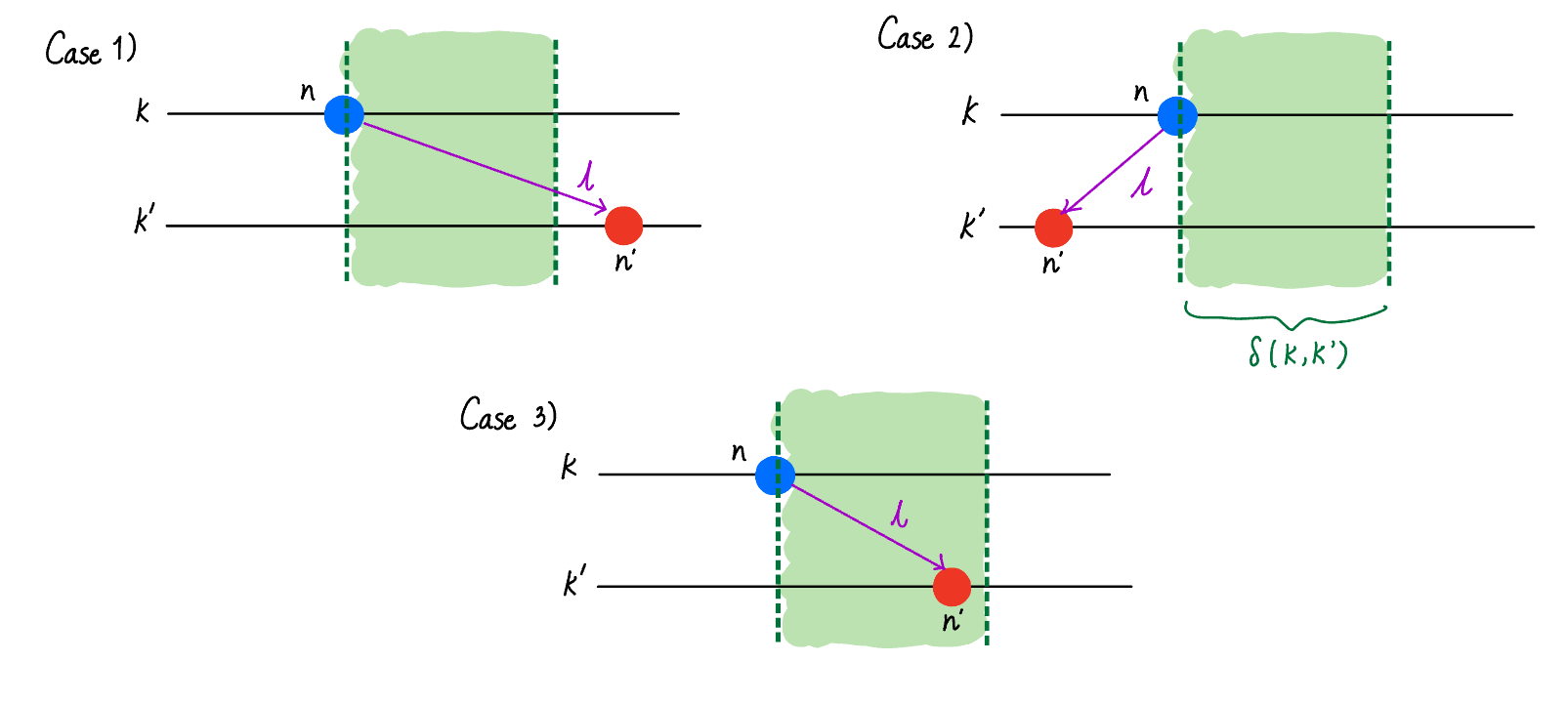}
    \caption{Finding the earliest available ground-departure node}
    \label{fig:cases_forA_L}
\end{figure}

While the reduction is mostly straightforward to apply, some care must be taken when dealing with wrap-around arcs.
In particular, consider the three cases in Figure~\ref{fig:cases_forA_L}.
This figure visualizes light travel from arrival-ground node $n \in N_R$ at terminal $k \in K$ to ground-departure node $n' \in N_E$ at terminal $k' \in K$.
The `earliest' ground-departure node is identified by starting at $n$, adding the travel time $\delta(k, k')$ and then moving forward in time until a ground-departure node is encountered at terminal $k$.
In Case 1, this results in a standard forward arc within the horizon.
In Case 2, the plan wraps around the horizon before a ground-departure node is encountered, and hence the resulting arc must be added to the set of wrap-around arcs $S$.
Finally, Case 3 is deceptive: because $n'$ cannot be reached within the travel time $\delta(k, k')$, this ground-departure node is only reached after the plan wraps around.
As such, this arc must also be added to $S$.
The same logic applies when identifying the `latest' node in the origin filtering step.

\subsection{Computation Performance Comparison}
\begin{figure}[H]
    \centering
    \includegraphics[width=0.8\linewidth]{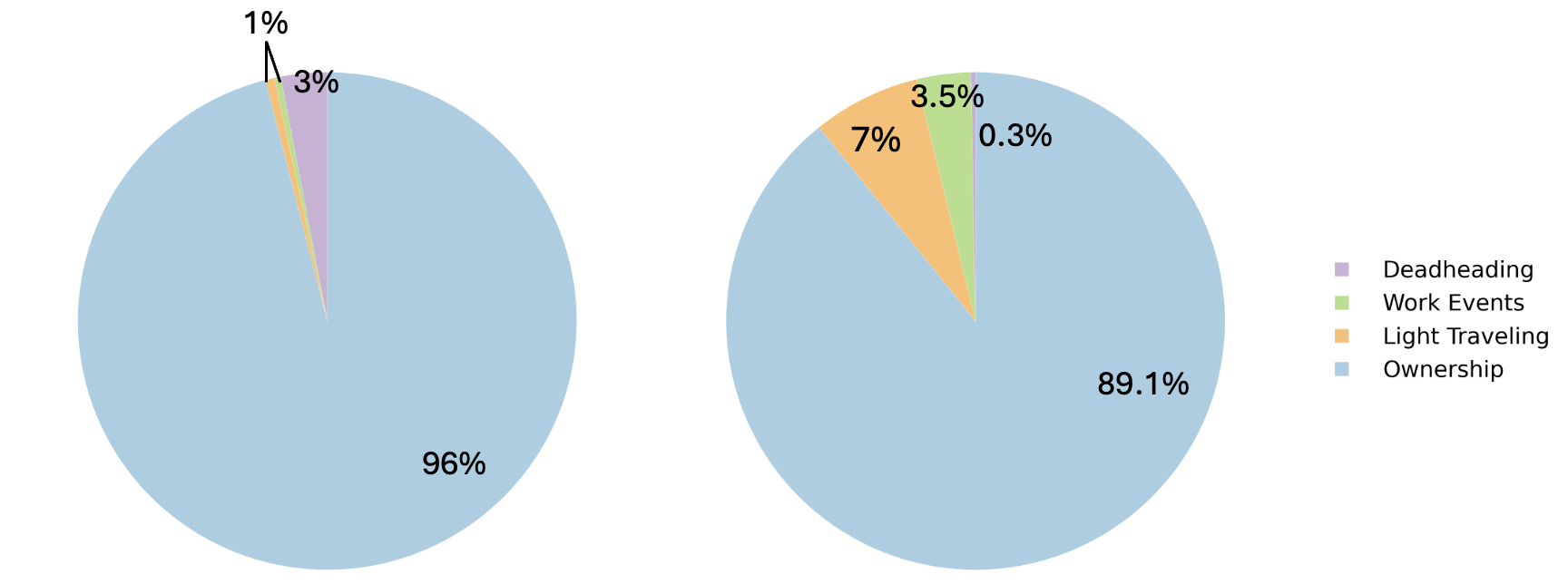}
    \caption{Cost Compositions in the Objective Value: MCF (Left) and Exact (Right)}
    \label{cost_composition}
\end{figure}

Figure \ref{cost_composition} shows objective cost breakdown for each method. In both, ownership cost dominates due to a high cost coefficient.  
The exact method consistently achieves a smaller fleet size with lower ownership cost.
This is offset by higher reliance on light traveling usage and work events, which facilitate more locomotive repositioning with a limited fleet. 
On the other hand, the MCF method solution exhibits a heavier use of deadheading as means of repositioning.

Table \ref{tab:time_composition} illustrates the share of weekly time spent on each activity. Each value is computed by dividing the total time locomotives spend on an activity (aggregated across the week) by the product of fleet size and the length of planning horizon. As illustrated, the smaller fleet size in the exact method leads to a higher proportion of active time, since total active time is fixed by the input data. 
The exact method also allocates a greater share of time to light travel, while MCF shows higher percentage of time spent in deadheading. 
Additionally, the exact model yields lower proportion of idle time, suggesting locomotives are less frequently sitting unused in inventory--consistent with the need to maintain system balance using fewer assets. This potentially suggests more efficient fleet utilization under a smaller inventory, a preferred solution depending on user preference.

\begin{table}[H]
    \centering
    \caption{Weekly Locomotive Activity Breakdown}
    \begin{tabular}{c c c c c c c c} \hline  
         Method&  Active&  DH&  Pre-Dept &  Post-Arvl&  Connection&  LT & Idle \\ \hline  
         Exact&  $53.9\%$&  $5.8\%$&  $8.2\%$&  $2.5\%$&  $14.8\%$ &  $1.7\%$ &  $13.2\%$ \\   
         MCF&  $51.9\%$&  $7.5\%$&  $8.2\%$&  $2.5\%$&  $15\%$&  $0.2\%$& $14.9\%$\\  \hline
    \end{tabular}
    \label{tab:time_composition}
\end{table}

\section{Trade-off Analysis} \label{Appendix:6.2}
 
\begin{figure}[H]
    \centering
    \subfigure[Varying Ownership Cost ($q$)]{
        \includegraphics[trim=0 0 270 0, clip, width=0.45\textwidth]{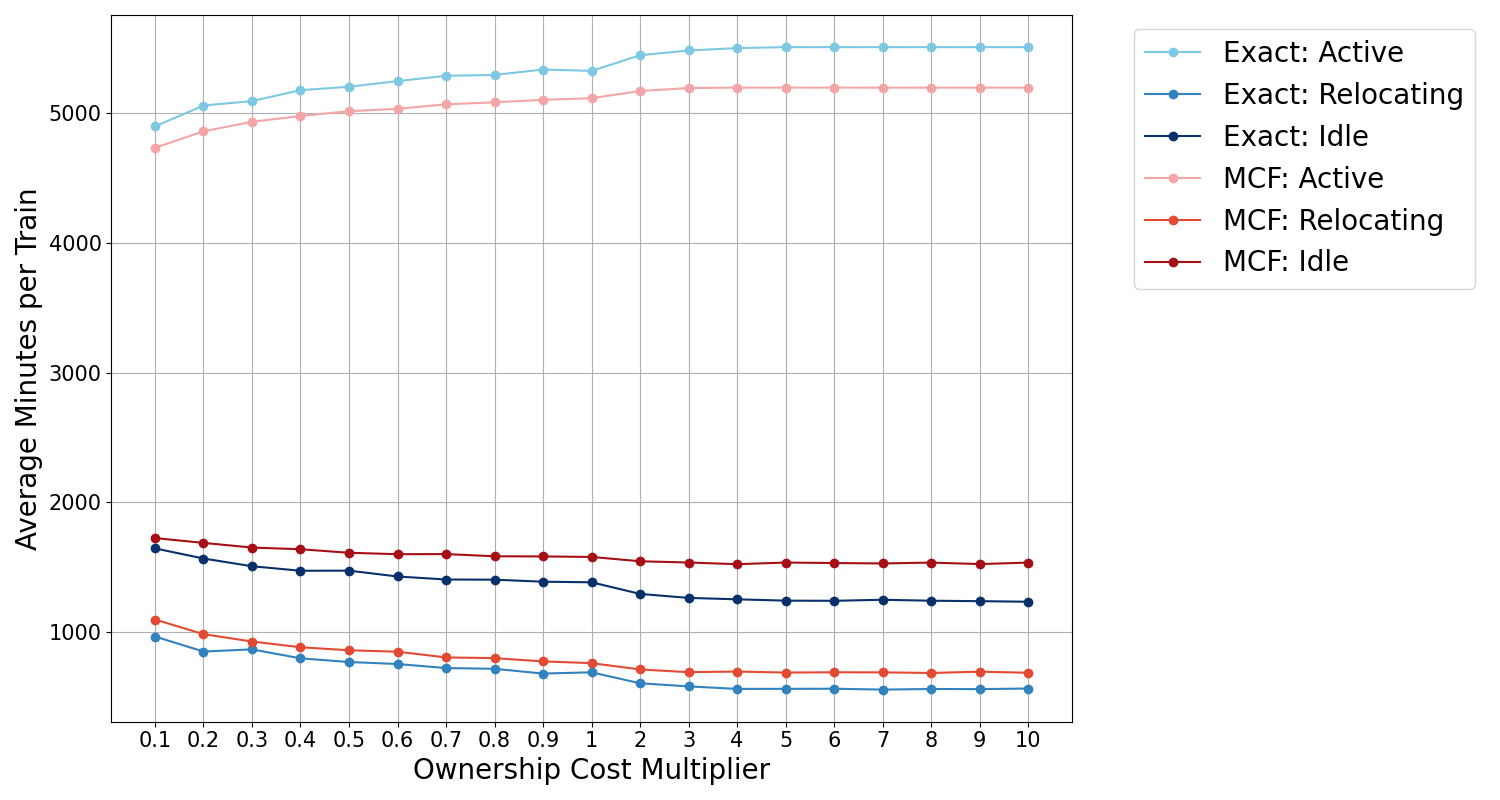}
        \label{fig:FC_util}
    }
    \hfill
    \subfigure[Varying Crew Cost ($e$)]{
        \includegraphics[trim=0 0 270 0, clip, width=0.45\textwidth]{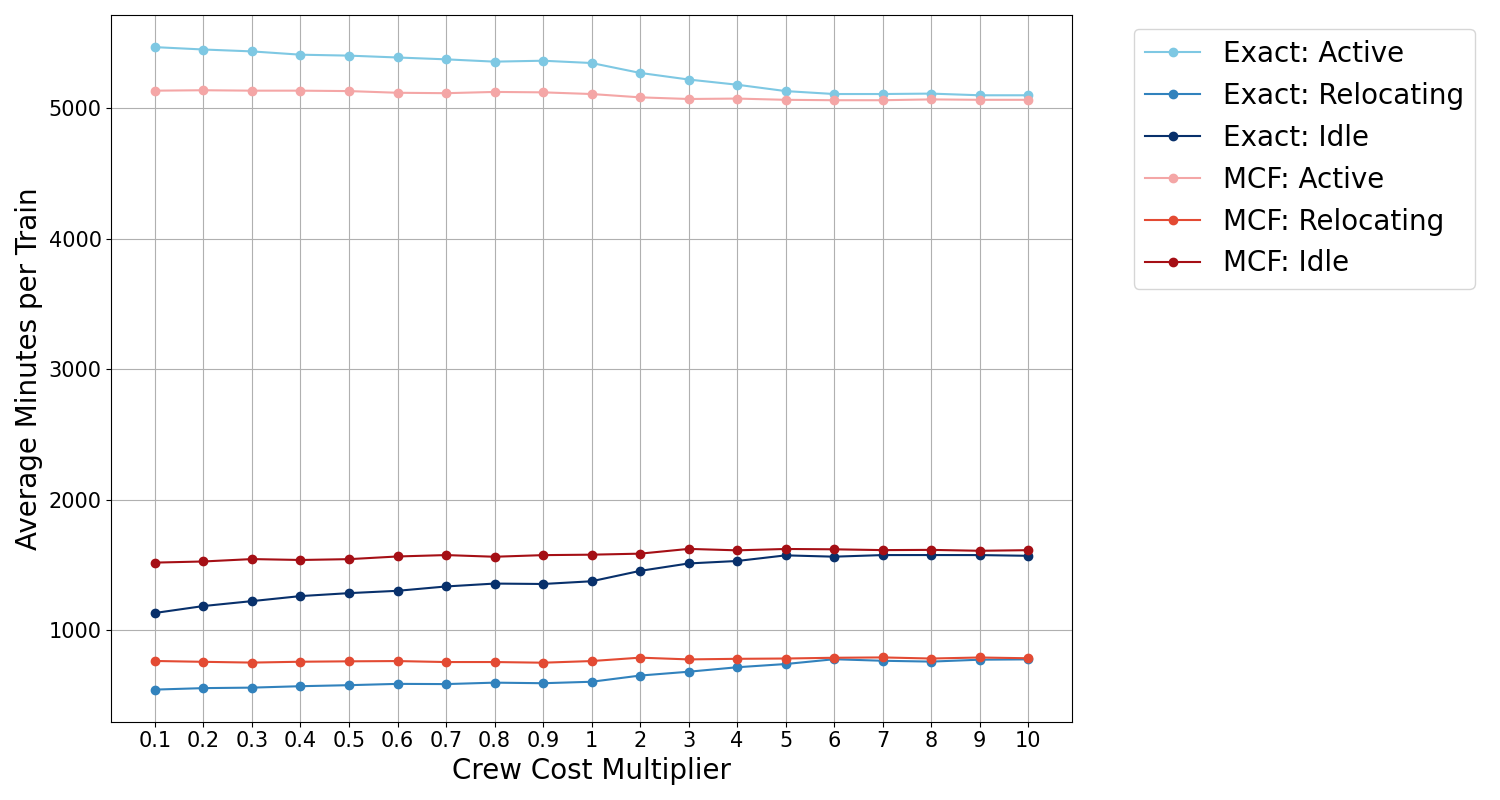}
        \label{fig:CC_util}
    }
    \hfill
    \subfigure[Varying Work Events Cost ($c$)]{
        \includegraphics[width=0.6\textwidth]{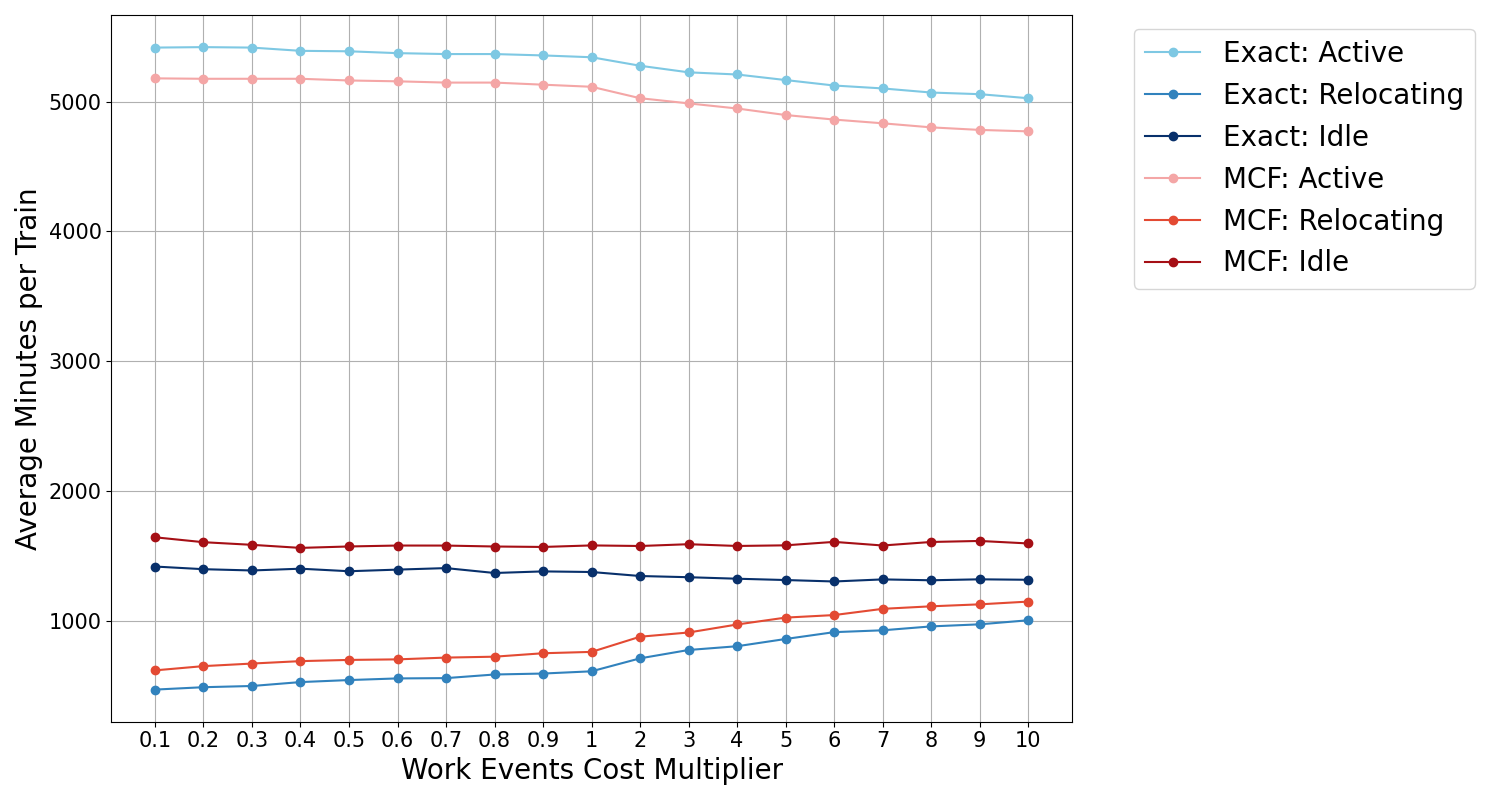}
        \label{fig:WC_util}
    }
    \caption{Weekly Average Minutes per Train Allocated to Locomotive Activities}
    \label{fig:utilization_avgTime}
\end{figure}

The trade-off between fleet sizing and locomotive utilization is further reflected in how locomotive time is allocated across activities. 
Figure \ref{fig:utilization_avgTime} summarizes average time spent in each activity under varying levels of $q$, $e$, and $c$. 
As fleet contracts (e.g., under high-$q$ scenarios), locomotives exhibit higher asset productivity, spending more time in active operation and less time idle.

Together, the trade-off analysis demonstrates how cost parameters shape optimization behavior and offers a spectrum of viable options. Rather than prescribing a single “best” plan, the models support informed decision-making by mapping out trade-offs between fleet size, light travel, and decentralization.
Crucially, these elements cannot be optimized in isolation. Their interactions are complex and nonlinear, and effective strategies must be designed at the system level. For decision-makers, this means jointly considering investments in fleet size, terminal flexibility, and repositioning policies. A clear understanding of these interdependencies is essential for crafting cost-effective and resilient operations.

\section{Model Extensions: Work Event Restrictions} \label{Appendix:6.3}
Table \ref{tab:alpha_summary} summarizes the versions and their control parameter used in the model extensions. The detailed model formulations are as follows.

\begin{table}[t]
\centering
\caption{Summary of Activation Control Parameters Across Model Variants}
\begin{tabular}{>{\raggedright}p{1.4cm} p{2cm} p{2.6cm} p{2cm} p{1cm} p{2.2cm} p{1.8cm}}
\toprule
\textbf{Variant} & \textbf{Parameter} & \textbf{Activation}& \textbf{Initial Value} & \textbf{Incre- ment}& \textbf{Constraints} & \textbf{Type} \\
\midrule
V1 & \(\lambda\) & N/A & 0 & +1 & (\ref{constr:verB_1})–(\ref{constr:verB_3}) & \shortstack{Capacity\\Increase} \\
V2 & \(\alpha^C\) & terminals \hspace{0.5cm} (\(k \in K^I\))& 0 & +1 & (\ref{constr:verCD_1})-(\ref{constr:verCD_2}), (\ref{constr:verC_1})–(\ref{constr:verC_3}) & \shortstack{Incremental \\ Expansion } \\
\addlinespace
V3 & \(\alpha^D\) & terminal-days (\((k,d) \in KD^I\))& 0 & +5 & (\ref{constr:verCD_1})-(\ref{constr:verCD_2}), (\ref{constr:verD_1})–(\ref{constr:verD_3}) & \shortstack{Incremental \\ Expansion }\\
\addlinespace
V4 &\(\alpha^E\) &  terminals \hspace{0.5cm} (\(k \in K\))& \(|K \setminus K^I|\) & +1& (\ref{constr:verE_1})–(\ref{constr:verE_3}) & \shortstack{Network \\ Redesign} \\
\addlinespace
V5 &\(\alpha^F\) & terminal-days (\((k,d) \in KD\))& \(|KD \setminus KD^I|\) & +5& (\ref{constr:verF_1})–(\ref{constr:verF_3}) & \shortstack{Network \\ Redesign} \\
\bottomrule
\end{tabular}
\label{tab:alpha_summary}
\end{table}

\subsection{Extension of Baseline Plan} 
\paragraph{V1: Increasing Capacity at Baseline-Active Terminal-Days}
V1 examines the impact of incrementally increasing capacity at terminal-day pairs already utilized under the current plan (${KD} \setminus {KD}^I$). A parameter $\lambda \in \mathbb{N}_0$ defines the maximum number of additional events allowed per such pair, and variants are expressed as $\text{V1}_\lambda$. This extension does not enforce adherence to any specific event allocations; rather, it expands an upper bound limit to work events.

The following constraints are included:
\begin{align}
    & \sum_{l \in \mathcal{L}_{kd}} y^{so}_l + y^{pu}_l \leq h_{k, d} + \lambda &\qquad \forall (k, d) \in KD \setminus {KD}^I \label{constr:verB_1}\\
    & \sum_{l \in \mathcal{L}_{kd}} y^{so}_l + y^{pu}_l \leq \theta &\qquad \forall (k, d) \in KD \setminus {KD}^I \label{constr:verB_2}\\
    & \sum_{l \in \mathcal{L}_{kd}} y^{so}_l + y^{pu}_l = 0 &\qquad \forall (k, d) \in {KD}^I \label{constr:verB_3}
\end{align}
Here, $\sum_{l \in \mathcal{L}_{kd}} y^{so}_l + y^{pu}_l$ captures the number of set-outs and pick-ups at terminal $k$ on day $d$.
Constraints (\ref{constr:verB_1}) permit up to $\lambda$ additional events at already active terminal-day pairs, while constraints (\ref{constr:verB_2}) enforce the global upper bound $\theta$. Constraints (\ref{constr:verB_3}) prohibit any activity at previously inactive terminal-day pairs. 
Varying $\lambda$ from 0 to $\theta - h_{k, d}$ gauges how much of the increased capacity will result in utilization, and thus, the maximum performance gain attainable from extending on the current plan.

\subsection{Incremental Expansion from Baseline Plan}
\paragraph{V2: Activating Baseline-Inactive Terminals}
This extension enables the activation of previously unused terminals, in addition to the active terminal-day pairs ($KD \setminus {KD}^I$). Hereafter, the baseline scenario denoted as V1' (when control parameter values are 0) is defined by constraints (\ref{constr:verCD_1})-(\ref{constr:verCD_2}). They restrict the capacity at active terminals to the minimum of twice the original allocation $h_{k,d}$ or the daily upper bound $\theta$. These depict a reasonable operational flexibility at active terminals, and represents the reference solution against which the following variants are evaluated.
\begin{align}
    & \sum_{l \in \mathcal{L}_{kd}} \left( y^{so}_l + y^{pu}_l \right) \leq 2 h_{k,d} && \forall (k, d) \in KD \setminus {KD}^I \label{constr:verCD_1} \\
    & \sum_{l \in \mathcal{L}_{kd}} \left( y^{so}_l + y^{pu}_l \right) \leq \theta && \forall (k, d) \in KD \setminus {KD}^I \label{constr:verCD_2} 
\end{align}

As $\alpha^C \in \mathbb{N}_0$ increases, previously inactive terminals are activated to extend the spatial scope of the baseline plan. Binary variables $z^1$ are defined for inactive terminals $k \in K^I$, indicating whether a terminal is made available. Constraints (\ref{constr:verC_1}) limit such activations to $\alpha^C$, while constraints (\ref{constr:verC_2}) ensure events occur only at activated terminals. This implies that as long as a terminal is activated, then it has ability to execute work events at any day of the week.
\begin{align}
    & \sum_{k \in K^I} z^1_k \leq \alpha ^C \label{constr:verC_1}\\
    & \sum_{l \in \mathcal{L}_{kd}} \left( y^{so}_l + y^{pu}_l \right) \leq \theta z^1_k && \forall (k, d) \in K^I \times D \label{constr:verC_2} \\
    & z^1_k \in \{ 0,1 \} && \forall k \in K^I \label{constr:verC_3}
\end{align}
This setup quantifies the incremental benefits gained from activating inactive terminals by gradually increasing $\alpha^C$. It studies how much improvements can be realized by expanding the active terminal set beyond the baseline configuration, offering insights into which terminals would provide the highest cost savings.

\paragraph{V3: Activating Baseline-Inactive Terminal-Day}
In this version, both spatial and temporal flexibility are integrated. The baseline scenario corresponds to aforementioned V1' and assumes no new terminal-day pair activations ($\alpha^D = 0$). When $\alpha^D > 0$, the model introduces additional terminal-day pair activations through binary variables $z^2$, defined over the set ${KD}^I$. Constraints (\ref{constr:verD_1}) allow up to $\alpha_D$ new activations, and constraints (\ref{constr:verD_2}) state work events assigned to an activated pair do exceed the upper bound $\theta$.
\begin{align}
    & \sum_{(k, d) \in {KD}^I} z^2_{kd} \leq \alpha ^D \label{constr:verD_1}\\
    & \sum_{l \in \mathcal{L}_{kd}} \left( y^{so}_l + y^{pu}_l \right) \leq \theta z^2_{kd} && \forall (k, d) \in {KD}^I  \label{constr:verD_2}\\
    & z^2_{kd} \in \{ 0,1 \} && \forall (k, d) \in {KD}^I \label{constr:verD_3}
\end{align}
Similar to V2, this approach navigates incremental adjustments to an established structure rather than overhauling an entire plan. It is distinct by optimizing the specific days of operations along with terminals, rather than full utilization of selective terminals. This diversifies the terminals that are activated, by offering narrowly targeted adjustments at the day-of-week level.

\subsection{Network-wide Redesign}
\paragraph{V4: Full Terminal Activation Flexibility} 
This formulation completely disregards the baseline plan and optimizes terminal activations across the full set $K$. A binary variable $\omega^1_k$ is introduced for each terminal $k \in K$, where $\omega^1_k = 1$ permits work events at terminal $k$. Constraints (\ref{constr:verE_1}) limit the number of activated terminals to $\alpha^E \in \mathbb{N}_0$. Constraints (\ref{constr:verE_2}) cap total events at each activated terminal to a maximum of $\theta$. Activated terminals are eligible to conduct work events on any day.  
\begin{align}
    & \sum_{k \in K} \omega^1_k \leq \alpha ^E \label{constr:verE_1}\\
    & \sum_{l \in \mathcal{L}_{kd}} \left( y^{so}_l + y^{pu}_l \right) \leq \theta \omega^1_k && \forall k \in K \label{constr:verE_2}\\
    & \omega^1_k \in \{ 0,1 \} && \forall k \in K \label{constr:verE_3}
\end{align}
This formulation constructs and identifies the most valuable terminals under an idealized scenario. It functions as a clean-slate benchmark, particularly relevant for network design. Comparing the resulting terminals to $K \setminus K^I$ reveals insights into the efficiency of the current plan.

\paragraph{V5: Full Terminal-Day Activation Flexibility} 
V5 is another clean-slate approach but incorporates decision-making at the terminal-day pair level. Binary variables $\omega^2_{k,d}$ indicate whether terminal $k$ can have work events on day $d$. Constraints (\ref{constr:verF_1}) bound the total number of active pairs by $\alpha^F$, and constraints (\ref{constr:verF_2}) limit events at each pair by $\theta$. 
\begin{align}
    & \sum_{(k, d) \in {KD}} \omega^2_{kd} \leq \alpha ^F \label{constr:verF_1}\\
    & \sum_{l \in \mathcal{L}_{kd}} \left( y^{so}_l + y^{pu}_l \right) \leq \theta \omega^2_{kd} && \forall (k, d) \in {KD}  \label{constr:verF_2}\\
    & \omega^2_{kd} \in \{ 0,1 \} && \forall (k, d) \in {KD} \label{constr:verF_3}
\end{align}
This approach offers the highest level of configurational flexibility. It accommodates a broad range of scenarios, from tightly constrained activation policies to near full coverage. It is particularly useful for testing robustness under disruptions such as reduced capacity (\emph{e.g.,} low $\alpha^E$ value), terminal congestion/closure, uncertain demand, or infrastructure disruptions. It can be a diagnostic tool to understand performance degradation under adversity, and for identifying critical thresholds in network resilience. Furthermore, it guides the design of minimal yet sufficient activation strategies to maintain feasibility and service continuity.

\subsection{Additional Computational Experiments}
\begin{figure}[H]
    \centering
    \includegraphics[width=0.55\linewidth]{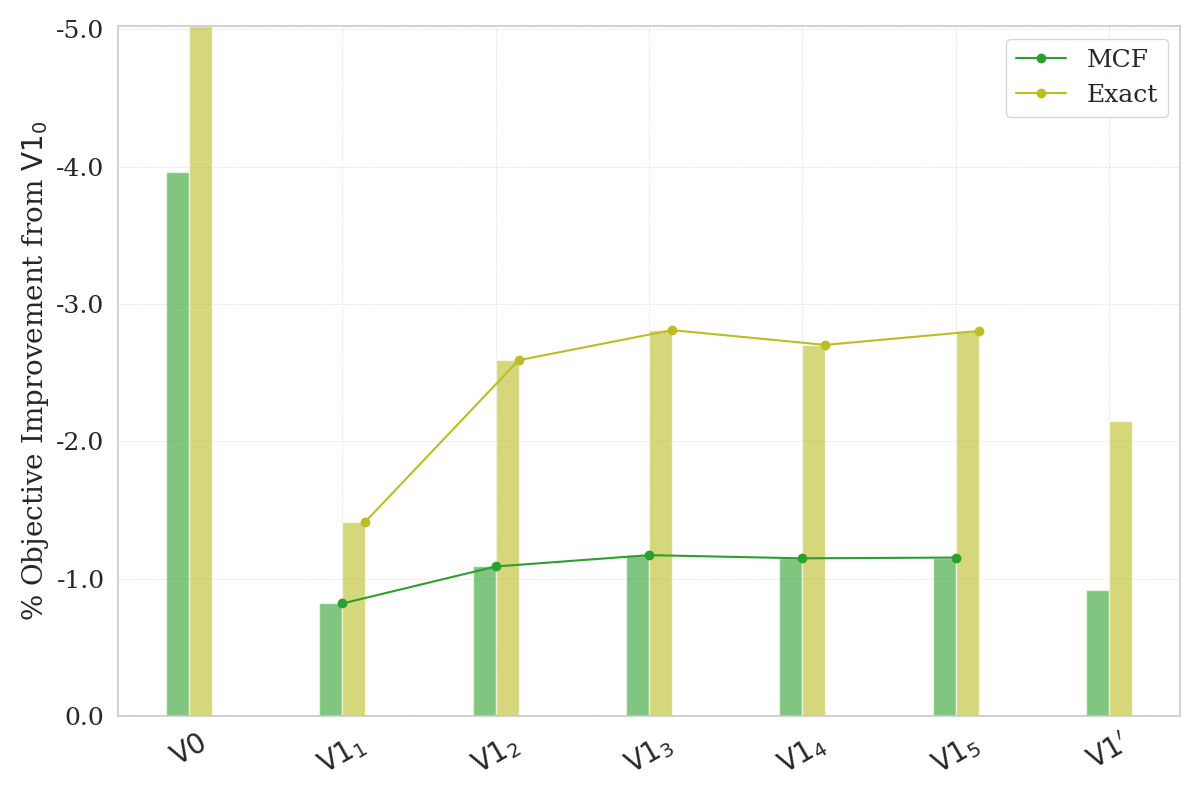}
    \caption{Objective Value Difference of Scenario V1 (\%) from $\text{V1}_0$}
    \label{fig:v1_obj_changes}
\end{figure}

$\text{V1}_0$ refers to the constrained model with no capacity increase ($\lambda=0$): it only reallocates work events within the baseline plan without increasing volume.
Figure \ref{fig:v1_obj_changes} reports the percent reduction in objective value relative to $\text{V1}_0$ as $\lambda$ increases. A one-unit increase cuts cost by ~$1.5\%$ and two units by ~$2.8\%$. Beyond that, marginal returns diminish, implying most optimization potential is captured with modest extension. Even a 1\% reduction translates to hundreds of thousands of dollars in weekly savings. 

\begin{figure}[h]
    \centering    \includegraphics[width=0.7\linewidth]{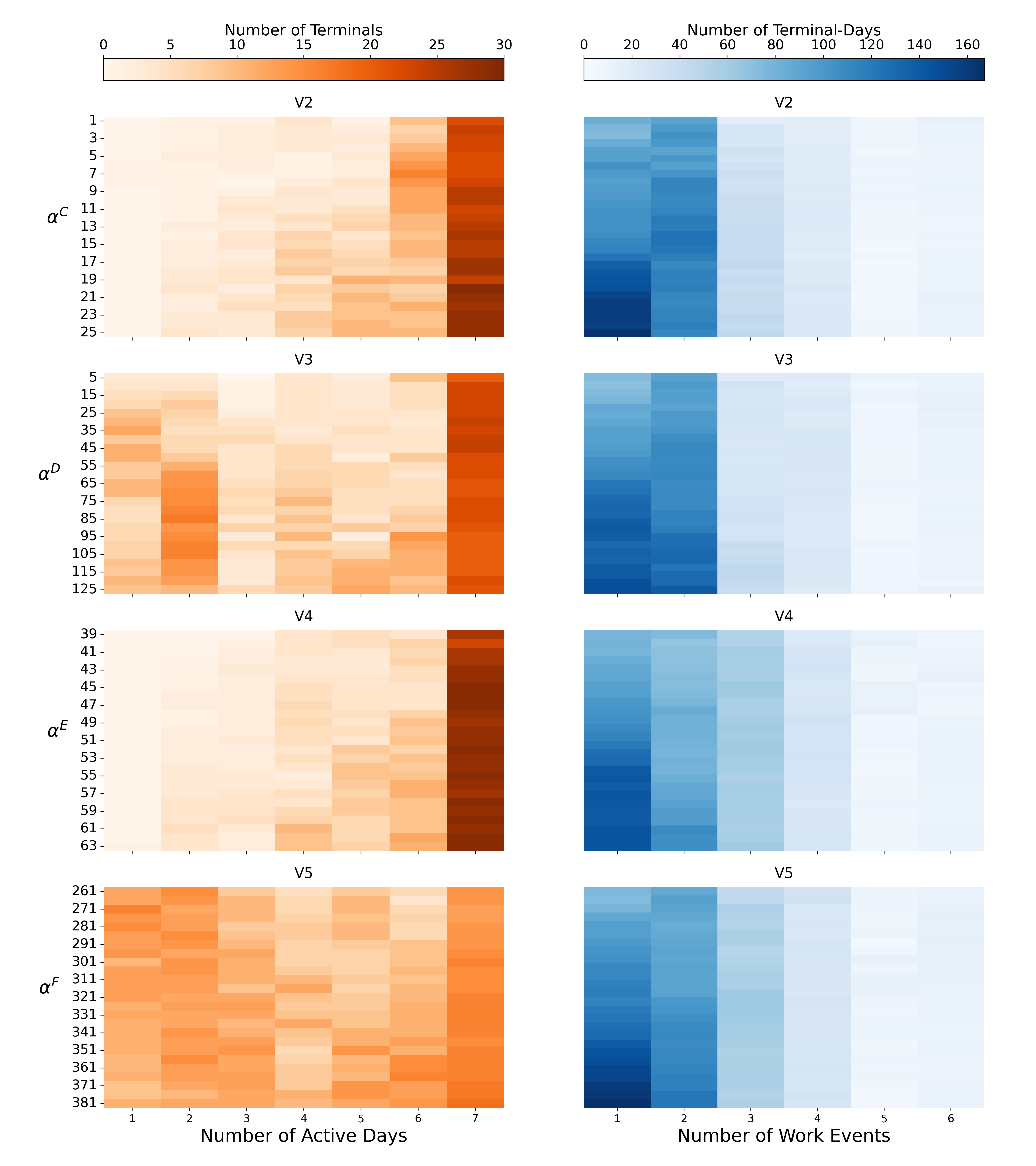}
    \caption{Spatio-Temporal Activation Patterns Across Model Variants Under MCF method}
    \label{fig:MCF_heatmaps}
\end{figure}

\subsection{Practical Implications}
Efficiency gains can be achieved without increasing the number of work events. A practical first step is to remove work events from baseline-active terminal-days that are consistently excluded from optimized plans across variants. This frees up resources with minimal disruption at no cost.

Incremental models (V2, V3) consistently remove a stable subset of baseline elements with fewer persisting as thresholds increase. Even retained terminals exhibit reduced day-level usage in V2, suggesting uneven marginal utility across days. V3 initially prunes a large share of baseline-active assignments that it takes several iterations to match the baseline's total active terminal-days count. 
Meanwhile, a core set of high-performing terminals (or terminal-days)--from both the baseline and new additions--recurs across all variants and thresholds.
These persistent selections form a resilient backbone for robust network function. Improvements come from dynamic assignments around this stable, high-impact core.

\begin{figure}[htbp]
    \centering
    \subfigure{
        \includegraphics[width=0.45\textwidth]{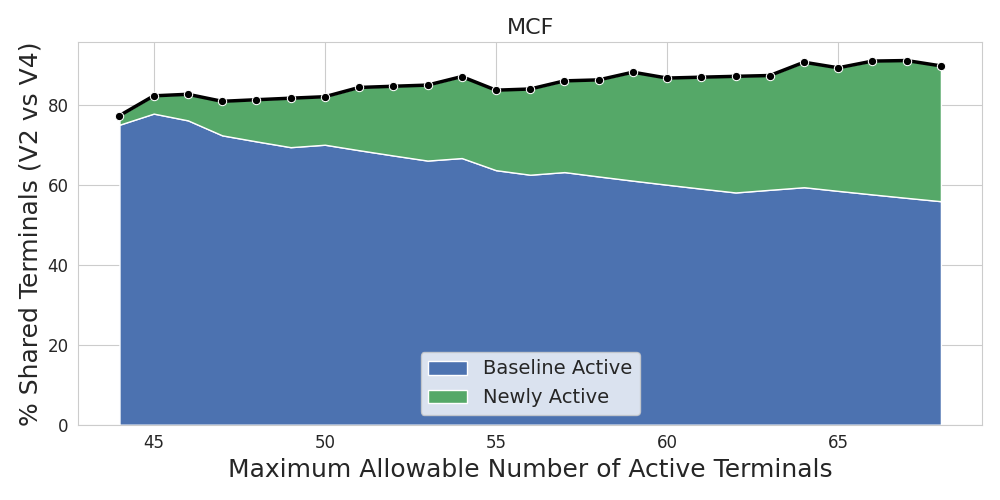}
        \label{fig:overlap_percentages_1}
    }
    \hfill
    \subfigure{
        \includegraphics[width=0.45\textwidth]{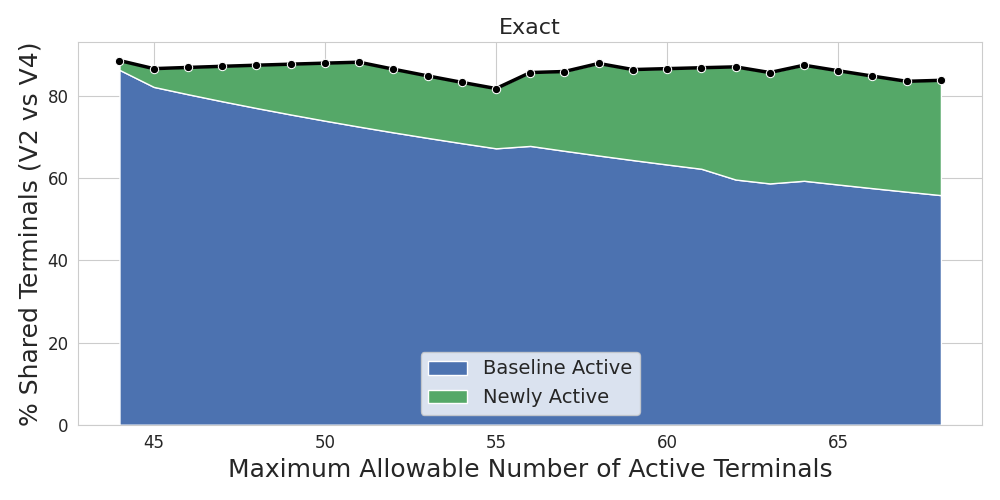}
        \label{fig:overlap_percentages_2}
    }

    \vspace{0.1cm}

    \subfigure{
        \includegraphics[width=0.45\textwidth]{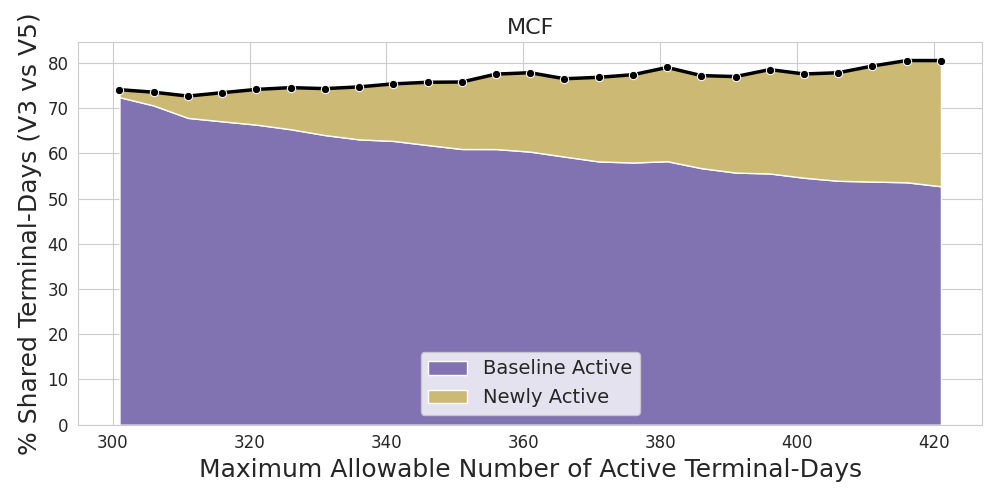}
        \label{fig:figoverlap_percentages_3}
    }
    \hfill
    \subfigure{
        \includegraphics[width=0.45\textwidth]{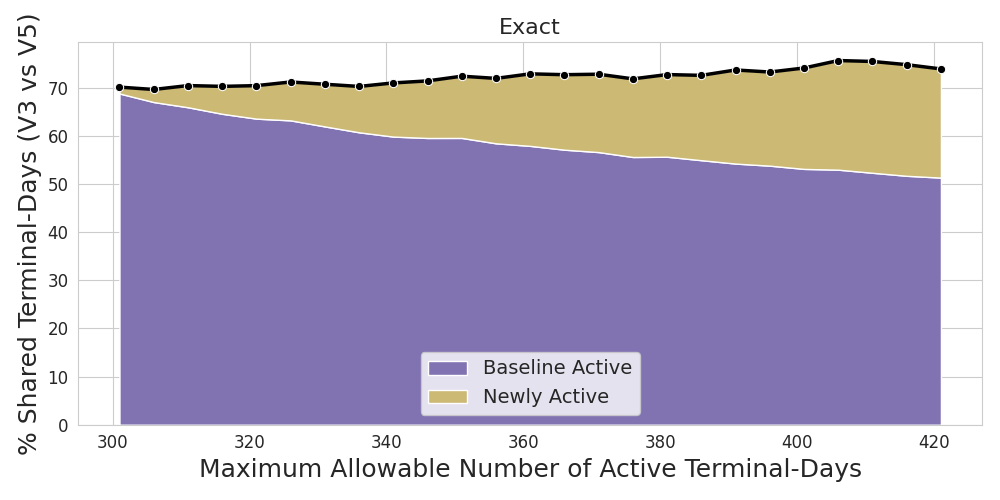}
        \label{fig:overlap_percentages_4}
    }

    \caption{Overlap in Selected Terminals and Terminal-Days Between Incremental and Redesign Models}
    \label{fig:overlap_percentages}
\end{figure}

Despite the theoretical flexibility of redesigns, roughly 80\% of V4 terminals and 70\% of V5 terminal-days overlap with their incremental counterparts (Figure \ref{fig:overlap_percentages}). 
This suggests clean-slate designs restructure selectively rather than complete overhaul. Performance gains are driven by the remaining 20-30\% that enable improved coordination with shared components. 
Also note that at lower activation levels, warm-starting constrains the model’s ability to meaningfully deviate from the baseline.

Lastly, effective work events generate cost savings not only through broader geographic reach but also through precisely timed execution. While relying on a few high-capacity terminals may seem appealing, results show that dynamically allocating work event capacity across both terminals and days consistently yields greater savings. Unlocking this potential may require operational flexibility through investments in staffing or infrastructure.

\end{document}